\documentclass[a4paper, 11pt]{article}

\usepackage[left=2.5cm, right=2.5cm, top=2cm, bottom=2.5cm]{geometry}
\usepackage[utf8]{inputenc} 
\usepackage[T1]{fontenc}    
\usepackage{hyperref}       
\usepackage{url}            
\usepackage{booktabs}       
\usepackage{amsmath, amssymb}     
\usepackage{nicefrac}       
\usepackage{microtype}      
\usepackage{lipsum}		
\usepackage{graphicx}
\usepackage{doi}
\usepackage{authblk}
\newtheorem{remark}{Remark}
\usepackage{enumerate}
\usepackage[section]{placeins}

\title{Inverse scattering of periodic surfaces with a superlens}

\date{May 5, 2023}	

\author[1]{Peijun Li \thanks{\href{mailto:lipeijun@math.purdue.edu}{lipeijun@math.purdue.edu}}}
\author[2,3]{Yuliang Wang \thanks{Corresponding author: \href{mailto:yuliangwang@bnu.edu.cn}{yuliangwang@bnu.edu.cn}}}
\affil[1]{Department of Mathematics, Purdue University, West Lafayette, Indiana 47907, USA}
\affil[2]{Research Center for Mathematics, Beijing Normal University, Zhuhai 519087, China}
\affil[3]{BNU-HKBU United International College, Zhuhai 519087, China}

\begin{document}

\maketitle

\begin{abstract} 
We propose a scheme for imaging periodic surfaces using a superlens. By employing an inverse scattering model and the transformed field expansion method, we derive an approximate reconstruction formula for the surface profile, assuming small amplitude. This formula suggests that unlimited resolution can be achieved for the linearized inverse problem with perfectly matched parameters. Our method requires only a single incident wave at a fixed frequency and can be efficiently implemented using fast Fourier transform. Through numerical experiments, we demonstrate that our method achieves resolution significantly surpassing the resolution limit for both smooth and non-smooth surface profiles with either perfect or marginally imperfect parameters.
\end{abstract}

\section{Introduction}

Resolution is a crucial factor in any wave imaging system. In traditional imaging systems, such as optical microscopes, the resolution is limited by approximately half of the wavelength, a rule known as the Rayleigh criterion or the Abbe diffraction limit. Decreasing the wavelength may not always be feasible, especially in source imaging problems, and high-frequency incident waves may cause side effects such as sample damage. Therefore, it is desirable to surpass the diffraction limit, enabling high resolution to be achieved with low-frequency waves.

Typically, the diffraction limit arises from the inability to capture evanescent waves, which convey fine details of the imaged object. However, these waves are confined near the object's surface, and their amplitude diminishes rapidly as the distance from the object increases. Consequently, a viable strategy for overcoming the diffraction limit involves recovering evanescent waves through measurements taken in close proximity to the object. This approach has been employed in near-field optical microscopy techniques, such as scanning near-field optical microscopy \cite{betzig1987collection,chen2019modern,hecht2000scanning,heinzelmann1994scanning} and photon scanning tunneling microscopy \cite{reddick1990photon}.

Optical microscopy techniques generate images of samples directly on the measuring device, often resulting in the loss of depth information or its recovery at a lower resolution. To reconstruct the full profile of a sample, it is desirable to solve inverse scattering problems based on the underlying wave equations. For infinite surfaces, the inverse scattering problem involves reconstructing the surface profile using scattered waves measured on a surface. One prevalent approach is to construct an objective functional derived from the discrepancy between measured and computed data for a given surface profile. The surface profile is then identified as a minimizer of the objective functional through optimization algorithms \cite{bao2013near,bao2016shape,zhang2017imaging,li2017kirsch}. While these methods offer versatility in terms of objective functionals and optimization algorithms, they are computationally expensive and prone to issues related to local minima.

An alternative approach, known as sampling methods, involves designing an indicator function with high contrast values across the surface profile. The surface can then be visually identified by plotting the indicator function within a sampling region \cite{ding2017imaging,li2021linear,hu2019non,liu2018direct}. While sampling methods offer the advantages of speed and minimal reliance on prior knowledge, they typically necessitate a substantial volume of measurement data, which may not be readily available in practical scenarios. Furthermore, these methods are intrinsically qualitative, leading to potential limitations in the quantitative precision of the resulting images.

In specific practical situations, surfaces may display a low profile with respect to period and wavelength. In such instances, techniques based on the transformed field expansion (TFE) can be utilized. The TFE was devised in \cite{nicholls2004josaa, nicholls2004improved} to address direct rough surface scattering problems. It was initially employed in \cite{baoli2013} to solve inverse scattering problems on infinite rough surfaces and has since been expanded to various contexts, encompassing periodic surfaces \cite{cheng2013near}, obstacles \cite{li2015nearObstacle,li2015numerical}, interior cavities \cite{li2015nearCavity}, penetrable surfaces \cite{bao2014near}, elastic waves \cite{li2015inverseElastic,li2015nearObstacleElastic}, and three-dimensional problems \cite{li2016nearElastic,bao2014inverse,jiang2017inverse}. These methods leverage the TFE and linearization to derive explicit and accessible reconstruction formulas. They are not only straightforward but also effective and efficient, necessitating merely a single incident field. Another advantage of these techniques is their capacity to attain convergence results and error estimates \cite{bao2014convergence,li2016convergence}, attributable to the explicit and accessible reconstruction formulas.

In this paper, we re-examine the near-field inverse scattering problem for periodic surfaces, aiming to further improve imaging resolution by integrating a superlens into the model. An ideal superlens is a slab of material characterized by a relative permittivity $\varepsilon = - 1$ and relative permeability $\mu = - 1$. The notion of materials with negative permittivity and permeability was initially proposed by Veselago in \cite{veselago1967electrodynamics}, who contended that such materials could exhibit atypical optical properties, including a negative refractive index. Subsequently, Pendry introduced the concept of employing a slab of this material as an imaging lens and illustrated the potential for eliminating the resolution limit, coining the term "superlens" in \cite{pendry2000negative}.

Following this, a multitude of numerical and experimental demonstrations of superlenses for electromagnetic waves have been showcased in the literature \cite{fang2005sub,moussa2005negative,aydin2007subwavelength}. In more recent developments, the concept has been expanded to include acoustic waves and elastic waves. Nevertheless, the general imaging approach in these studies remains akin to that of a traditional microscope, where light from a source passes through the lens, ultimately forming a two-dimensional projection of the imaged profile on the image plane. In contrast, our imaging scheme relies on inverse scattering, providing a comprehensive reconstruction of the entire profile.

In this study, we propose an imaging model that comprises an impenetrable periodic surface, a slab of negative index material situated above the surface, a measurement surface atop the slab, and a single incident field. By employing the low-profile assumption and the TFE method, we derive an explicit reconstruction formula. We illustrate that, in the ideal case, unlimited resolution can be achieved, while high resolution can be obtained with slightly non-ideal parameters. To the best of our knowledge, this research represents the first instance of integrating negative index material into inverse scattering problems, thereby unlocking the full-profile imaging resolution enhancement potential offered by the superlens.

This work can be viewed as an extension of our previous study \cite{bao2016nearFar}, wherein we attained enhanced resolution utilizing a slab with a high refractive index. The enhancement factor in \cite{bao2016nearFar} is approximately proportional to the slab's refractive index; however, it tends toward infinity in this paper as the parameters approach ideal values. Furthermore, we apply the TFE differently in this work by treating the entire domain as a single boundary value problem. As a result, the reconstruction formula becomes substantially simpler, facilitating preliminary resolution analysis. Moreover, we commence with the most general settings for the slab's material parameters to derive a comprehensive formula, which encompasses the no-slab, high-index slab, and negative-index slab as special cases.

The remainder of the paper is organized as follows. In Section \ref{sec:model}, we establish the physical model and derive the boundary value problem for the entire domain using interface conditions and a transparent boundary condition. Section \ref{sec:TFE} presents the TFE, which reduces the problem into a recursive system of one-dimensional boundary value problems. In Section \ref{sec:analytical}, we obtain the analytical solutions for the leading and linear terms. By solving the linearized inverse problem, we derive the reconstruction formula and conduct elementary resolution analysis in Section \ref{sec:reconstruction}. Numerical experiments are presented in Section \ref{sec:numerical} to demonstrate the effectiveness of the proposed imaging scheme. Finally, Section \ref{sec:conclusion} offers a conclusion and discusses future research directions.

\section{The model problem }\label{sec:model}

Consider the electromagnetic scattering by a periodic surface defined as 
\begin{equation*}
\Gamma_f = \{ (x, y) \in \mathbb{R}^2: y = f (x) \},
\end{equation*}
where $f$ is a periodic function with period $\Lambda$. We assume that $f$ can be expressed as
\begin{equation}
f (x) = \delta g (x), \label{eq:49}
\end{equation}
where $0 < \delta \ll 1$ is a small constant, referred to as the surface deformation parameter. The surface $\Gamma_f$ is assumed to be perfectly electrically conducting (PEC), and the half space above it is assumed to be a vacuum.

\begin{figure}[h]
\centering
  \includegraphics[width=0.35\textwidth]{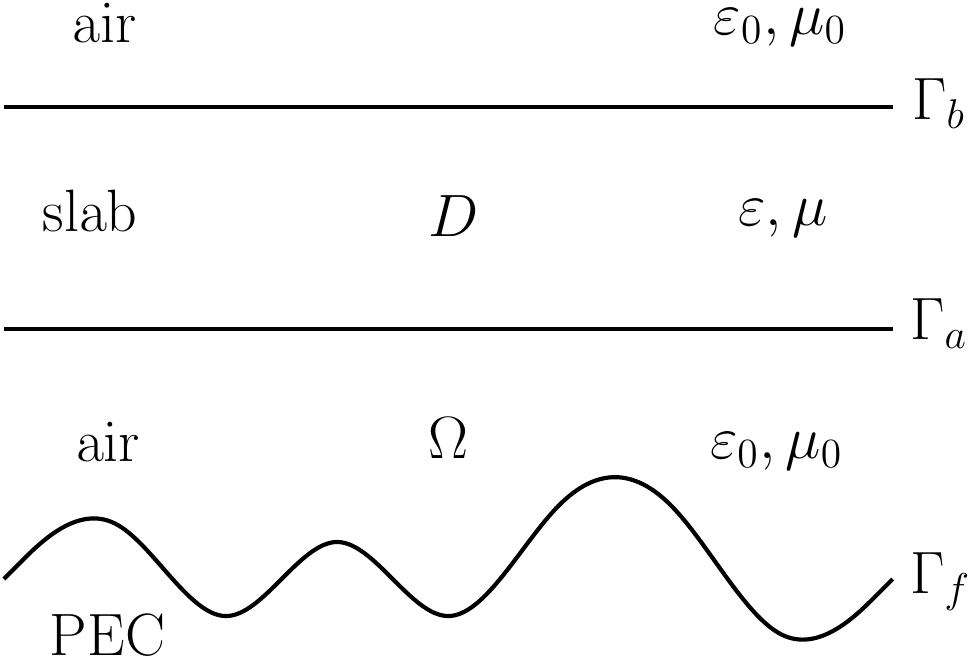}
  \caption{Geometry of the model problem. The slab is positioned with its lower boundary $\Gamma_h$ and upper boundary $\Gamma_b$ above the scattering surface $\Gamma_f$, and the measurements are taken on $\Gamma_b$.}
  \label{fig:geometry}
\end{figure}

We introduce a slab of double negative index material above the scattering surface, with its two surfaces defined by $\Gamma_a = \{ (x, a): x \in \mathbb{R} \}$ and $\Gamma_b = \{ (x, b): x \in \mathbb{R} \}$. The domain bounded between $\Gamma_f$ and $\Gamma_a$ is denoted by $\Omega = \{ (x, y) \in \mathbb{R}^2: f (x) < y < a \}$, and the domain between $\Gamma_a$ and $\Gamma_b$ is represented by $D = \{ (x, y) \in \mathbb{R}^2: a < y < b \}$. The problem geometry is illustrated in Fig. \ref{fig:geometry}.

The propagation of electromagnetic waves is described by the time-harmonic Maxwell's equations: 
\begin{equation}
\label{eq:12} \nabla \times \mathbf{H} + i \omega \varepsilon \mathbf{E} =
0, \quad \nabla \times \mathbf{E} - i \omega \mu \mathbf{H} = 0,
\end{equation}
where $\mathbf{E}$ and $\mathbf{H}$ are the electric and magnetic fields, respectively, $\omega$ represents the angular frequency, $\varepsilon$ and $\mu$ are the permittivity and permeability, respectively. In this work, we consider the transverse electric (TE) mode, assuming all variables are independent of the $z$ coordinate, and $\mathbf{E} = [0, 0,
u]^\top$, $\mathbf{H} = [H_1, H_2, 0]^\top$. Then, Eq. \eqref{eq:12} reduces to
\begin{equation}
\label{eq:19} (\Delta + \omega^2 \varepsilon \mu) u = 0.
\end{equation}

The relative permittivity and permeability are denoted by $\varepsilon_r = \varepsilon / \varepsilon_0$ and $\mu_r = \mu / \mu_0$, respectively, where $\varepsilon_0$ and $\mu_0$ are the permittivity and permeability of vacuum, respectively. Define the free-space wavenumber $\kappa = \omega \sqrt{\varepsilon_0 \mu_0}$, then Eq. \eqref{eq:19} can be rewritten as
\begin{equation*}
  (\Delta + \eta^2) u = 0,
\end{equation*}
where $\eta = \kappa \sqrt{\varepsilon_r \mu_r}$. The positive $x$-axis is taken as the branch cut for any square root in this paper. The term $\sqrt{\varepsilon_r \mu_r}$ is commonly referred to as the refractive index.

A plane wave with an incident angle $\theta \in (- \pi / 2, \pi / 2)$ is described by $u_{\mathrm{in}} (x, y) = e^{i (\alpha x - \beta y)}$, where $\alpha = \kappa \sin \theta$ and $\beta = \kappa \cos \theta$. However, we will focus on normal incidence, i.e., $\theta = 0$. In this case, the incident field simplifies to
\begin{equation}
u_{\mathrm{in}} (x, y) = e^{- i \kappa y}, \label{eq:56}
\end{equation}
where both $u_{\mathrm{in}}$ and $u$ are periodic in $x$ with a period of $\Lambda$. 

Given a scattering surface and incident field, the direct scattering problem aims to determine the total field $u(x, y)$. The near-field imaging problem we consider in this paper is an inverse scattering problem, seeking to determine the scattering surface function $f(x)$ from the measured data of the total field at $\Gamma_b$.

The incident field given by Eq. \eqref{eq:56} clearly satisfies the Helmholtz equation:
\begin{equation*}
  (\Delta + \kappa^2) u_{\mathrm{in}} (x, y)
= 0, \quad y > b.
\end{equation*}
Let $u_{\mathrm{sc}} = u - u_{\mathrm{in}}$ denote the scattered field, which satisfies the same Helmholtz equation:
\begin{equation}
\label{eq:37} (\Delta + \kappa^2) u_{\mathrm{sc}} (x, y)
= 0, \quad y > b.
\end{equation}
Assuming that $u_{\mathrm{sc}}$ consists only of upward propagating and evanescent waves, we can obtain the Rayleigh expansion from Eq. \eqref{eq:37} as follows:
\begin{equation}
u_{\mathrm{sc}} = \sum_{n \in \mathbb{Z}}
u_{\mathrm{sc}}^{(n)} (b) e^{i [\alpha_n x + \beta_n (y-b)]}, \quad y \geq b, \label{eq:57}
\end{equation}
where $u_{\mathrm{sc}}^{(n)} (b)$ denotes the Fourier coefficient of $u_{\mathrm{sc}} (x, b)$, and
\begin{equation}
\alpha_n = \frac{2 \pi n}{\Lambda}, \quad \beta_n = \sqrt{\kappa^2 -
\alpha_n^2}. \label{eq:59}
\end{equation}
Taking $\partial_y$ in Eq. \eqref{eq:57} and evaluating at $\Gamma_b$, we get
\begin{equation}
\partial^+_y u_{\mathrm{sc}}
=\mathcal{T}u_{\mathrm{sc}}^+ \quad \text{on}~\Gamma_b,
\label{eq:60}
\end{equation}
where the boundary operator $\mathcal{T}$ is defined by
\begin{equation*}
\mathcal{T}v = \sum_{n \in \mathbb{Z}} i \beta_n v^{(n)},
\end{equation*}
and the "$+$" sign indicates taking partial derivative or limit from above. It is straightforward to verify
\begin{equation}
\partial_y^+ u_{\text{in}} =\mathcal{T}u_{\text{in}}^+ + \rho \quad
\text{on}~\Gamma_{b}, \label{eq:63}
\end{equation}
where
\begin{equation}
\rho = - 2 i \kappa e^{- i \kappa b}. \label{eq:58}
\end{equation}
Adding Eq. \eqref{eq:63} to Eq. \eqref{eq:60} yields the transparent boundary condition:
\begin{equation}
\partial_y^+ u =\mathcal{T}u^+ + \rho \quad \text{on}~ \Gamma_b. 
\label{eq:64}
\end{equation}

From the interface conditions $\nu\times\mathbf E^+=\nu\times\mathbf E^-$ and $\nu\times\mathbf{H}^+=\nu\times\mathbf{H}^-$ on $\Gamma_b$ with $\nu=[0, 0, 1]^\top$, we obtain
\begin{equation}
\label{eq:33} u^- = u^+, \quad \frac{1}{\mu_r} \partial_y^- u = \partial_y^+
u \quad \text{on}~ \Gamma_b,
\end{equation}
where the "$-$" sign indicates taking partial derivative or limit from below. Combining Eqs. \eqref{eq:64} and \eqref{eq:33} leads to 
\begin{equation*}
  \frac{1}{\mu_r} \partial_y^- u =\mathcal{T}u + \rho \quad
\text{on}~ \Gamma_b.
\end{equation*}
Similarly, we have
\begin{equation*}
  u^+ = u^-, \quad \frac{1}{\mu_r} \partial_y^+ u = \partial_y^-
u \quad \text{on}~ \Gamma_a.
\end{equation*}

For simplicity of notations, we will denote the relative permittivity and permeability by $\epsilon$ and $\mu$, respectively, in the following. To summarize, we obtain the following boundary-interface value problem for the total field $u$ in $\overline{D \cup \Omega}$:
\begin{subequations}
\label{eq:bivp-original}
\begin{alignat}{2}
& \partial_y^- u = \mu (\mathcal{T} u + \rho) & \quad &
\text{on}~ \Gamma_b, \\
& (\Delta + \eta^2) u = 0 & \quad &
\text{in}~ D, \\
& u^+ = u^-, \quad \partial_y^+ u = \mu \partial_y^- u &
\quad & \text{on}~ \Gamma_a, \\
& (\Delta + \kappa^2) u = 0 & \quad & \text{in}~ \Omega,\\
& u = 0 & \quad & \text{on}~ \Gamma_f. 
\end{alignat}
\end{subequations}

\section{Transformed field expansion}\label{sec:TFE}

We apply the TFE method, along with power series and Fourier series expansions, to reduce the governing equations \eqref{eq:bivp-original} to a recursive system of ordinary differential equations. This approach has been previously employed in \cite{baoli2013} for the purpose of near-field imaging of infinite rough surfaces.

Consider the change of variables
\begin{align*}
  \tilde{x} & = x, \quad \tilde{y} = a \left( \frac{y - f}{a - f} \right),
  \quad (x, y) \in \Omega,\\
  \tilde{x} & = x, \quad \tilde{y} = y, \quad (x, y) \in D,
\end{align*}
which transforms the domain $\Omega$ to the rectangle $\Omega_0 = (0, \Lambda)
\times (0, a)$, the surface $\Gamma_a$ to itself, the surface $\Gamma_f$ to
the plane $\Gamma_0=\{(x, 0): x\in\mathbb R\}$, and the domain $D$ to itself. It is important to note that the change of variables is defined in the whole computational domain, as opposed to focusing solely on the domain $\Omega$ as in \cite{bao2016nearFar}. This approach allows us to derive significantly simpler solution forms and reconstruction formulas.

Applying the chain
rule leads to the following differentiation rules in $\Omega$:
\begin{align*}
  \partial_x &= \partial_{\tilde{x}} - f'  \left( \frac{a - \tilde{y}}{a - f}
  \right) \partial_{\tilde{y}},\\
  \partial_y &= \left( \frac{a}{a - f} \right) \partial_{\tilde{y}},\\
  \partial_{xx} &= \partial_{\tilde{x}  \tilde{x}} + (f')^2  \left( \frac{a -
  \tilde{y}}{a - f} \right)^2 \partial_{\tilde{y}  \tilde{y}} - 2 f'  \left(
  \frac{a - \tilde{y}}{a - f} \right) \partial_{\tilde{x}  \tilde{y}} \\
  &- \left[
  f''  \left( \frac{a - \tilde{y}}{a - f} \right) + 2 (f')^2  \frac{(a -
  \tilde{y})}{(a - f)^2} \right] \partial_{\tilde{y}},\\
  \partial_{yy} &= \left( \frac{a}{a - f} \right)^2 \partial_{\tilde{y} 
  \tilde{y}}.
\end{align*}

Let
\[ \tilde{u} (\tilde{x}, \tilde{y}) = u (x, y), \quad (x, y) \in \Omega \cup
   D. \]
Substituting the differentiation rules to the boundary-interface value problem
Eqs. \eqref{eq:bivp-original} and dropping the tilde over the variables for simplicity of notations, we obtain 
the transformed problem:
\begin{subequations}
\label{eq:bivp-transformed}
\begin{alignat}{2}
& \partial_y^{-} u=\mu( \mathcal{T} u+\rho) &\quad& \text{on}~ \Gamma_b, \\
& \left(\Delta+\eta^2\right) u=0 &\quad& \text {in}~ D,  \\
& u^{+}=u^{-}, \quad \left(1-\frac{f}{a}\right) \partial_y^{+} u=\mu \partial_y^{-} u &\quad& \text {on}~ \Gamma_a, \\
& \left(c_1 \partial_{x x}+c_2 \partial_{y y}+c_3 \partial_{x y}+c_4 \partial_y+c_1 \kappa^2\right) u=0 &\quad& \text {in}~ \Omega_0, \\
& u=0 &\quad& \text {on}~ \Gamma_0,
\end{alignat}
\end{subequations}
where
\begin{equation*}
    \begin{cases}
      c_1=(a-f)^2, \quad c_2=\left[(a-y) f^{\prime}\right]^2+a^2, \\
      c_3=-2(a-y)(a-f) f^{\prime}, \quad c_4=-(a-y)\left[(a-f) f^{\prime \prime}+2\left(f^{\prime}\right)^2\right].
    \end{cases} 
\end{equation*}

Based on the smallness assumption \eqref{eq:49}, we consider
the power series expansion
\begin{equation}
  u (x, y) = \sum_{m = 0}^{\infty} \delta^m u_m (x, y), \quad (x, y) \in
  \Omega_0 \cup D. \label{eq:15}
\end{equation}
Substituting Eqs. \eqref{eq:49} and \eqref{eq:15} into Eqs. \eqref{eq:bivp-transformed}
yields the recursive system of equations
\begin{subequations}
\label{eq:ibvp-recursive}
    \begin{alignat}{2}
        & \partial_y^{-} u_m=\mu\left( T u_m+\rho_m\right) &\quad& \text{on}~ \Gamma_b, \\
        & \left(\Delta+\eta^2\right) u_m=0 &\quad& \text{in}~ D,  \\
        & u_m^{+}=u_m^{-}, \quad \partial_y^{+} u_m=\mu \partial_y^{-} u_m+\tau_m &\quad& \text{on}~ \Gamma_a, \\
        & \left(\Delta+\kappa^2\right) u_m=v_m &\quad& \text{in}~ \Omega_0, \\
        & u_m=0 &\quad& \text{on}~ \Gamma_0,
    \end{alignat}
\end{subequations}
where
\begin{equation}
  v_m = D_1 u_{m - 1} + D_2 u_{m - 2}, \quad \tau_m = T u_{m - 1}^+, \quad
  \rho_m = \delta_{0,m} \rho, \label{eq:22}
\end{equation}
where $\delta_{i,j}$ is the Kronecker delta. The differential operators $D_1, D_2, T$ are given by
\begin{align*}
  D_1 & = \frac{1}{a}  [2 g \partial_{xx} + 2 g' (a - y) \partial_{xy} + g''
  (a - y) \partial_y + 2 \kappa^2 g],\\
  D_2 & = - \frac{1}{a^2}  \{ g^2 \partial_{xx} + (g')^2 (a - y)^2
  \partial_{yy} + 2 gg' (a - y) \partial_{xy} \\
  &- [2 (g')^2 - gg''] (a - y),
  \partial_y + \kappa^2 g^2 \},\\
  T & = \frac{g}{a} \partial_y,
\end{align*}
and $\rho$ is given by Eq. \eqref{eq:58}.

As the variables $u_m, v_m, \tau_m$, and $\rho_m$ are periodic in the $x$ direction with a period of $\Lambda$, the Fourier series expansion can be applied to the governing equations \eqref{eq:ibvp-recursive} to obtain a recursive system of boundary-interface value problems:
\begin{subequations}
\label{eq:ibvp-Fourier}
    \begin{alignat}{2}
        & \partial_y^{-} u_m^{(n)}=\mu\left[i \beta_n u_m^{(n)}+\rho_m^{(n)}\right], & y=b, \label{eq:30} \\
        & \partial_{y y} u_m^{(n)}+\gamma_n^2 u_m^{(n)}=0, & a<y<b, \label{eq:26}\\
        & \left(u_m^{+}\right)^{(n)}=\left(u_m^{-}\right)^{(n)}, \quad \partial_y^{+} u_m^{(n)}=\mu \partial_y^{-} u_m^{(n)}+\tau_m^{(n)}, & y=a, \label{eq:29}\\
        & \partial_{y y} u_m^{(n)}+\beta_n^2 u_m^{(n)}=v_m^{(n)}, & 0<y<a, \label{eq:25}\\
        & u_m^{(n)}=0, & y=0. \label{eq:27}
    \end{alignat}
\end{subequations}
where $\beta_n$ are given in Eq. \eqref{eq:59} and
\begin{equation}
  \gamma_n = \sqrt{\eta^2 - \alpha_n^2} . \label{eq:70}
\end{equation}

\section{Analytical solutions}\label{sec:analytical}

Solving Eqs. \eqref{eq:ibvp-Fourier} by the variation of parameters leads to the
general solution
\begin{equation}
  u_m^{(n)} (y) = \left\{\begin{array}{ll}
    a_{m, n} e^{i \gamma_n y} + b_{m, n} e^{- i \gamma_n y}, & \quad a < y <
    b,\\
    c_{m, n} e^{i \beta_n y} + d_{m, n} e^{- i \beta_n y} + \int_0^y \Phi_n
    (y, z) v_m^{(n)} (z) d z, & \quad 0 < y < a,
  \end{array}\right. \label{eq:31}
\end{equation}
where $a_{m, n}, b_{m, n}, c_{m, n}, d_{m, n}$ are coefficients to be
determined, and the integral kernel
\begin{equation*}
  \Phi_n (y, z) = \frac{\sin \beta_n  (y - z)}{\beta_n} \quad \text{if}~
  \beta_n \neq 0.
\end{equation*}
Applying the boundary and interface conditions \eqref{eq:30}, \eqref{eq:29},
and \eqref{eq:27}, we obtain the linear system
\begin{equation}
  \begin{pmatrix}
    i (\gamma_n - \mu \beta_n) e^{i \gamma_n b} & - i (\gamma_n + \mu \beta_n)
    e^{- i  \gamma_n b} & 0 & 0\\
    e^{i \gamma_n a} & e^{- i \gamma_n a} & - e^{i \beta_n a} & - e^{- i
    \beta_n a}\\
    i \gamma_n e^{i \gamma_n a} & - i \gamma_n e^{- i \gamma_n a} & - i \mu
    \beta_n e^{i \beta_n a} & i \mu \beta_n e^{- i \beta_n a}\\
    0 & 0 & 1 & 1
  \end{pmatrix} 
  \begin{pmatrix}
    a_{m, n}\\
    b_{m, n}\\
    c_{m, n}\\
    d_{m, n}
  \end{pmatrix} = 
  \begin{pmatrix}
    r_{m, n}\\
    s_{m, n}\\
    t_{m, n}\\
    0
  \end{pmatrix}, \label{eq:34}
\end{equation}
where
\begin{equation}
  r_{m, n} = \mu \rho_m^{(n)}, \quad s_{m, n} = \int_0^a \Phi_n (a, z)
  v_m^{(n)} (z) d z, \quad t_{m, n} = \mu \int_0^a \partial_y \Phi_n (a, z)
  v_m^{(n)} (z) d z + \tau_m^{(n)} . \label{eq:4}
\end{equation}
A direct calculation shows that the determinant of the matrix in Eq. \eqref{eq:34} is
given by
\begin{align}
  \phi_n & = (\gamma_n - \mu \beta_n) e^{i \gamma_n  (b - a)}  [(\gamma_n +
  \mu \beta_n) e^{i \beta_n a} - (\gamma_n - \mu \beta_n) e^{- i \beta_n a}]
  \nonumber\\
  & \quad + (\gamma_n + \mu \beta_n) e^{i \gamma_n  (a - b)}  [(\gamma_n + \mu
  \beta_n) e^{- i \beta_n a} - (\gamma_n - \mu \beta_n) e^{i \beta_n a}] .
  \label{eq:72} 
\end{align}

\subsection{Zeroth order term}

For $m = 0$ it follows from Eq. \eqref{eq:22} that
\[ 
  v_0 = 0, \quad \tau_0 = 0, \quad \rho_0^{(n)} = \rho \delta_{0, n}.
\]
The Fourier coefficients are given by
\[ r_{0, n} = \mu \rho \delta_{0, n}, \quad s_{0, n} = 0, \quad t_{0, n} = 0.
\]
Solving Eq. \eqref{eq:34} by Cramer's rule yields the coefficients
\begin{subequations}
\label{eq:coeffs-zero}
    \begin{align}
          a_{0, n} & = \frac{2 \mu \rho}{i \phi_0} e^{- i \eta a}  (\mu \kappa \cos
          \kappa a + i \eta \sin \kappa a) \delta_{0, n},  \\
          b_{0, n} & = - \frac{2 \mu \rho}{i \phi_0} e^{i \eta a}  (\mu \kappa \cos
          \kappa a - i \eta \sin \kappa a) \delta_{0, n},  \\
          c_{0, n} & = \frac{2 \eta \mu \rho}{i \phi_0} \delta_{0, n}, \quad d_{0, n}
          = - \frac{2 \eta \mu \rho}{i \phi_0} \delta_{0, n}.  
\end{align}
\end{subequations}
Substituting Eqs. \eqref{eq:coeffs-zero} into Eq. \eqref{eq:31}, we obtain the
leading term of the solution
\begin{equation}
  u_0 (x, y) = \left\{\begin{array}{ll}
    a_{0, 0} e^{i \eta y} + b_{0, 0} e^{- i \eta y}, & \quad a \leq y \leq
    b,\\
    \dfrac{4 \eta \mu \rho}{\phi_0} \sin \kappa y, & \quad 0 \leq y \leq a.
  \end{array}\right. \label{eq:39}
\end{equation}
It should be noted that $u_0$ is independent of $x$ and corresponds to the total field when the scattering surface is flat, i.e., when the surface height deviation $\delta$ is zero.

\subsection{First order term}

For $m = 1$ we have $\rho_1 = 0$ and hence $r_{1, n} = 0$. Solving
Eq. \eqref{eq:34} by Cramer's rule for $a_{1, n}, b_{1, n}$ yields
\begin{align*}
  a_{1, n} & = \frac{2 (\gamma_n + \mu \beta_n)}{\phi_n} e^{- i \gamma_n b} 
  [\mu \beta_n \cos (\beta_n a) s_{1, n} - \sin (\beta_n a) t_{1, n}],\\
  b_{1, n} & = \frac{2 (\gamma_n - \mu \beta_n)}{\phi_n} e^{i \gamma_n b} 
  [\mu \beta_n \cos (\beta_n a) s_{1, n} - \sin (\beta_n a) t_{1, n}] .
\end{align*}
Substituting Eq. \eqref{eq:4} with $m = 1$ into the above equations and using the
difference formula for the sine function, we obtain
\begin{equation}
  a_{1, n} = - \frac{2 (\gamma_n + \mu \beta_n)}{\phi_n} e^{- i \gamma_n b}
  \psi_n, \quad b_{1, n} = - \frac{2 (\gamma_n - \mu \beta_n)}{\phi_n} e^{i
  \gamma_n b} \psi_n, \label{eq:23}
\end{equation}
where
\begin{equation}
  \psi_n = \mu \int_0^a \sin (\beta_n y) v_1^{(n)} (y) d y + \sin (\beta_n a)
  \tau_1^{(n)} . \label{eq:47}
\end{equation}

By substituting Eq. \eqref{eq:39} into Eq. \eqref{eq:22} and taking the Fourier transform, we obtain the following expressions for $v_1^{(n)}(y)$ and $\tau_1^{(n)}$:
\begin{subequations}
\label{eq:v1-tau1}
    \begin{align}
        v_1^{(n)} (y) & = \frac{4 \kappa \eta \mu \rho}{a \phi_0} [2 \kappa \sin(\kappa y) - \alpha_n^2 (a - y) \cos (\kappa y)] g^{(n)}, \quad 0 < y < a, \\
    \tau_1^{(n)} & = \frac{4 \kappa \eta \mu^2 \rho}{a \phi_0} \cos (\kappa a) g^{(n)}.
\end{align}
\end{subequations}
By substituting Eqs. \eqref{eq:v1-tau1} into Eq. \eqref{eq:47} and using the identity $\kappa^2 = \alpha_n^2 + \beta_n^2$, we derive the expression for $\psi_n$ through a lengthy but straightforward calculation:
\begin{equation}
\psi_n = \frac{4 \kappa \eta \mu^2 \rho \beta_n}{\phi_0} g^{(n)} . \label{eq:48}
\end{equation}
Substituting Eq. \eqref{eq:48} into Eq. \eqref{eq:23} yields the following expressions for $a_{1, n}$ and $b_{1, n}$:
\begin{subequations}
    \label{eq:51}
    \begin{align}
        a_{1, n} &= - \frac{8 \kappa \eta \mu^2 \rho \beta_n}{\phi_0 \phi_n} (\gamma_n + \mu \beta_n) e^{- i \gamma_n b} g^{(n)}, \\
        b_{1, n} &= - \frac{8 \kappa \eta \mu^2 \rho \beta_n}{\phi_0 \phi_n} (\gamma_n - \mu \beta_n) e^{i \gamma_n b} g^{(n)}.
    \end{align}
\end{subequations}
Finally, substituting Eqs. \eqref{eq:51} into Eq. \eqref{eq:31} provides the Fourier coefficients of the linear term:
\begin{equation}
u_1^{(n)} (y) = - \frac{8 \kappa \eta \mu^2 \rho \beta_n}{\phi_0 \phi_n} [(\gamma_n + \mu \beta_n) e^{i \gamma_n (y - b)} + (\gamma_n - \mu \beta_n) e^{i \gamma_n (b - y)}] g^{(n)}, \quad a < y < b. \label{eq:52}
\end{equation}
Evaluating Eq. \eqref{eq:52} at $y = b$ results in the fundamental identity of this paper:
\begin{equation}
\label{eq:53} 
g^{(n)} = \Upsilon_n u_1^{(n)} (b),
\end{equation}
where the scaling factor is given by:
\begin{equation}
\label{eq:scaling} \Upsilon_n = - \frac{\phi_0 \phi_n}{16 \kappa \eta \mu^2 \rho \beta_n \gamma_n}.
\end{equation}

\begin{remark}
For the sake of conciseness, we omit specific cases where either $\beta_n = 0$ or $\gamma_n = 0$ from this paper. However, it is important to note that these cases can be addressed separately if necessary.
\end{remark}

\section{Reconstruction Formula}\label{sec:reconstruction}

By retaining only the zeroth and first order terms in the power series expansion given by Eq. \eqref{eq:15}, we arrive at the approximation:
\begin{equation}
u(x, y) \approx u_0(x, y) + \delta u_1(x, y). \label{eq:71}
\end{equation}
Evaluating Eq. \eqref{eq:71} at $y = b$ and applying the Fourier transform with respect to $x$ gives
\begin{equation}
u^{(n)}(b) \approx u_0^{(n)}(b) + \delta u_1^{(n)}(b). \label{eq:8}
\end{equation}
Substituting Eqs. \eqref{eq:53} and \eqref{eq:49} into Eq. \eqref{eq:8}, we deduce
\begin{equation}
f^{(n)} = \Upsilon_n [u^{(n)}(b) - u_0^{(n)}(b)]. \label{eq:65}
\end{equation}
This implies that the Fourier coefficients of the scattering surface function $f(x)$ can be computed directly from the Fourier coefficients of the total field $u(x, b)$. The scattering surface function is then reconstructed as follows:
\begin{equation*}
f(x) = \sum_{|n| \leq N} f^{(n)} e^{i \alpha_n x},
\end{equation*}
where $N \in \mathbb{Z}$ represents the cut-off frequency. 

By examining Eq. \eqref{eq:65}, we observe that the coefficients $\Upsilon_n$ play a significant role in characterizing the stability and resolution of the inverse scattering problem. These coefficients typically have complex forms. In the following section, we explore some special cases. Referring to the Abbe diffraction limit, we define the minimum feature size $d_n$ for the $n$-th Fourier mode of $f$ as half of the period of $e^{i \alpha_n x}$, which is expressed as $d_n = \Lambda / (2|n|)$. Additionally, we denote the wavelength of the incident field in free space as $\lambda = 2 \pi / \kappa$.

\subsection{No-slab Case}

First, we consider the case in which the slab is absent, corresponding to the parameter values $\epsilon = 1$ and $\mu = 1$. It follows that $\eta = \kappa$ and $\gamma_n = \beta_n$. Substituting these parameters into Eq. \eqref{eq:72}, we obtain
\begin{equation}
\phi_n = - 4 \beta_n^2 e^{- i \beta_n b} . \label{eq:73}
\end{equation}
Plugging Eq. \eqref{eq:73} into Eq. \eqref{eq:scaling} yields
\begin{equation*}
| \Upsilon_n | = \frac{1}{2 \kappa} |e^{- i \beta_n b} |,
\end{equation*}
which is consistent with the result presented in \cite{cheng2013near}. We distinguish two cases:

\begin{enumerate}
\item If $| \alpha_n | < \kappa$, i.e., $|n| < \Lambda / \lambda$, then $\beta_n$ is a real number, and $| \Upsilon_n | = 1 / (2 \kappa)$ remains constant. Therefore, we can achieve stable reconstruction of the corresponding Fourier modes $f^{(n)}$. Moreover, since the minimum feature size $d_n > \lambda / 2$, the resolution of these modes falls within the limits set by the Abbe diffraction limit.
\item If $| \alpha_n | > \kappa$, i.e., $|n| > \Lambda / \lambda$, then $\mathrm{Re}~\beta_n = 0$, $\mathrm{Im}~\beta_n = \sqrt{\alpha_n^2 - \kappa^2}$, and
\begin{equation*}
| \Upsilon_n | = \frac{1}{2 \kappa} e^{| \beta_n | b}
\end{equation*}
is approximately an exponential function of $|n|$. Thus, the reconstruction of the corresponding Fourier modes $f^{(n)}$ becomes increasingly unstable as $|n|$ increases.
\end{enumerate}
A comprehensive and precise analysis of errors in this scenario can be found in \cite{bao2014convergence}.

\subsection{Ordinary Dielectric Material}

Let us assume that the slab is made of an ordinary dielectric material with parameter values $\epsilon > 1$ and $\mu = 1$. This case was previously analyzed by \cite{bao2016nearFar}. The novel reconstruction formula presented in this paper, Eq. \eqref{eq:65}, establishes a one-to-one correspondence between the Fourier coefficients. This is in contrast to the convolutional type found in \cite{bao2016nearFar}. As a result, the formula is significantly simpler and more convenient to analysis.

\subsection{Double Negative Metamaterial}\label{sec:perfect}

In this section, we focus on the main topic of our paper, which is the examination of slabs that have negative permittivity $\varepsilon$ (or negative real part) and negative permeability $\mu$.

\subsubsection*{Perfectly Matched Parameters}

We start by assuming $\epsilon = \mu = - 1$, referred to as perfectly matched parameters, in line with the original concept of the superlens proposed by \cite{pendry2000negative}. As a result, we have $\eta = \kappa$ and $\gamma_n = \beta_n$. By substituting these parameters into Eq. \eqref{eq:72}, we obtain
\begin{equation}
  \phi_n = 4 \beta_n^2 e^{i \beta_n (b - 2 a)} . \label{eq:79}
\end{equation}
Substituting Eq. \eqref{eq:79} into Eq. \eqref{eq:scaling} gives
\begin{equation*}
  | \Upsilon_n | = \frac{1}{2 \kappa} | e^{i \beta_n (b - 2 a)} | .
\end{equation*}
\begin{enumerate}[(i)]
  \item If $b \geq 2 a$, then $| \Upsilon_n | \leq 1 / (2 \kappa)$ for all $n
\in \mathbb{Z}$. Therefore, we can reconstruct all frequency modes with stability, leading to unlimited resolution!
  \item If $b < 2 a$, then
  \[ | \Upsilon_n | = \frac{1}{2 \kappa} e^{| \beta_n |  (2 a - b)}. 
  \]
 Consequently, this scenario is equivalent to the case without a slab but with a measurement distance of $2 a - b$. Furthermore, since $2 a - b < b$, we can expect an improvement in resolution in this situation.
\end{enumerate}

\subsubsection*{Imperfect Parameters}

Ideal parameters are often unattainable in practical applications due to physical and engineering limitations. In particular, the lens constituent material is typically lossy, which leads to an effective permittivity that is a complex-valued quantity. To address this issue, we introduce the parameters 
\[
\epsilon = -1+\epsilon', \quad \mu =-1 + \mu',
\]
where $\varepsilon'$ and $\mu'$ are small-amplitude, real or complex numbers. Using Equation \eqref{eq:70}, we obtain $\gamma_n^2 = \beta_n^2 + \iota$, where $\iota = \kappa^2 (\epsilon' \mu' - \epsilon' - \mu')$. Let $n$ be a fixed integer from $\mathbb{Z}$. As $\iota \rightarrow 0$, we can easily confirm that $| \Upsilon_n | \to \dfrac{1}{2 \kappa} | e^{i \beta_n (b - 2 a)} |$, which is the scaling factor for perfectly matched parameters. This implies that the resolution limit approaches that of the perfectly matched parameters as $\varepsilon' \rightarrow 0$ and $\mu' \rightarrow 0$.

\section{Numerical Experiments}\label{sec:numerical}

\begin{figure}[t]
\centering
  \includegraphics[width=0.33\textwidth]{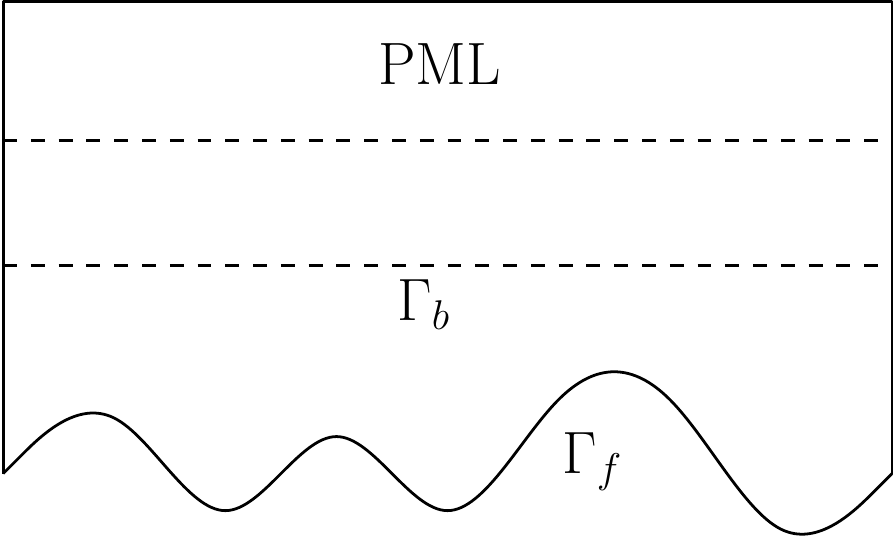}
  \caption{Sketch of the computational domain for the finite element method. A
  PML is placed above \(\Gamma_b\) to absorb the scattered field, allowing \(u = 0\)
  to be used as the boundary condition on the outer boundary of the PML.
  Periodic boundary conditions are enforced on the left and right boundaries, while 
  the usual PEC boundary condition \(u = 0\) is imposed on \(\Gamma_f\).
  \label{fig:FEM}}
\end{figure}

In the following numerical experiments, we keep the period fixed at \(\Lambda = 1\) and the wavelength at \(\lambda = 1.1\). To obtain the total field $u$ for a given scattering surface $f$, we solve the corresponding direct scattering problem using the finite element method along with the perfectly matched layer (PML) technique to truncate the infinite computational domain to a finite one. Fig. \ref{fig:FEM} shows a sketch of the computational domain.

The field data is then sampled at points uniformly distributed on \(\Gamma_b\), 
and random values are added to the samples to simulate measurement noise. Specifically,
the measurement data is given by 
\[ U = \left\{ u (x_m, b) (1 + r_m) : x_m = \frac{m}{M}, 0 \leq m \leq M
   \right\}, \]
where \(M = 100\), and \(r_m\) are independent random numbers drawn from the
uniform distribution in \([- 0.05, 0.05]\), i.e., the relative noise level is \(5\%\).
Throughout the experiment, the lower boundary of the slab is fixed at \(a = 0.1\) while the upper boundary is fixed at \(b = 0.2\).

\subsection{Smooth profile function}

The scattering surface function is defined as \(f (x) = \delta g (x)\), where
\(\delta = 0.01\), and the function
\begin{equation}
  g (x) = 0.4 \cos (2 \pi x) + 0.3 \cos (6 \pi x) + 0.2 \cos (20 \pi x)
  \label{eq:smoothProfile}
\end{equation}
represents frequency modes with \(|n| = 1, 3, 10\). The clear cut-off of the Fourier modes is used to test the proposed method's capabilities.

\begin{figure}[t]
  \begin{tabular}{ccc}
    \includegraphics[width=0.3\textwidth]{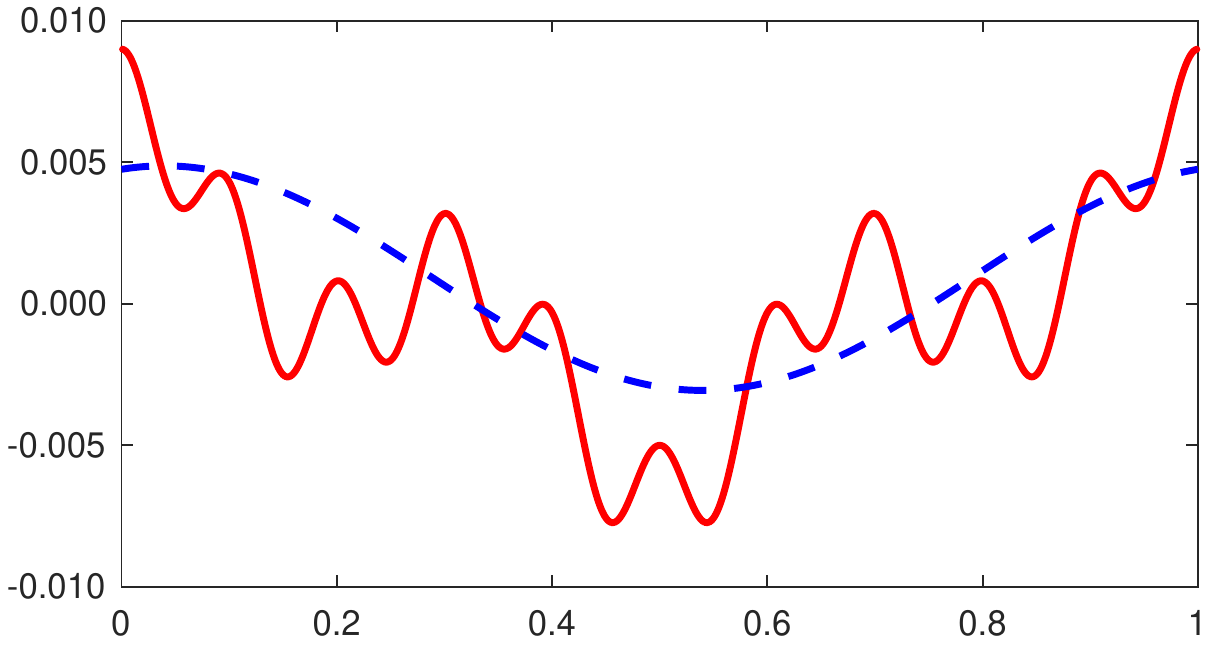} &
    \includegraphics[width=0.3\textwidth]{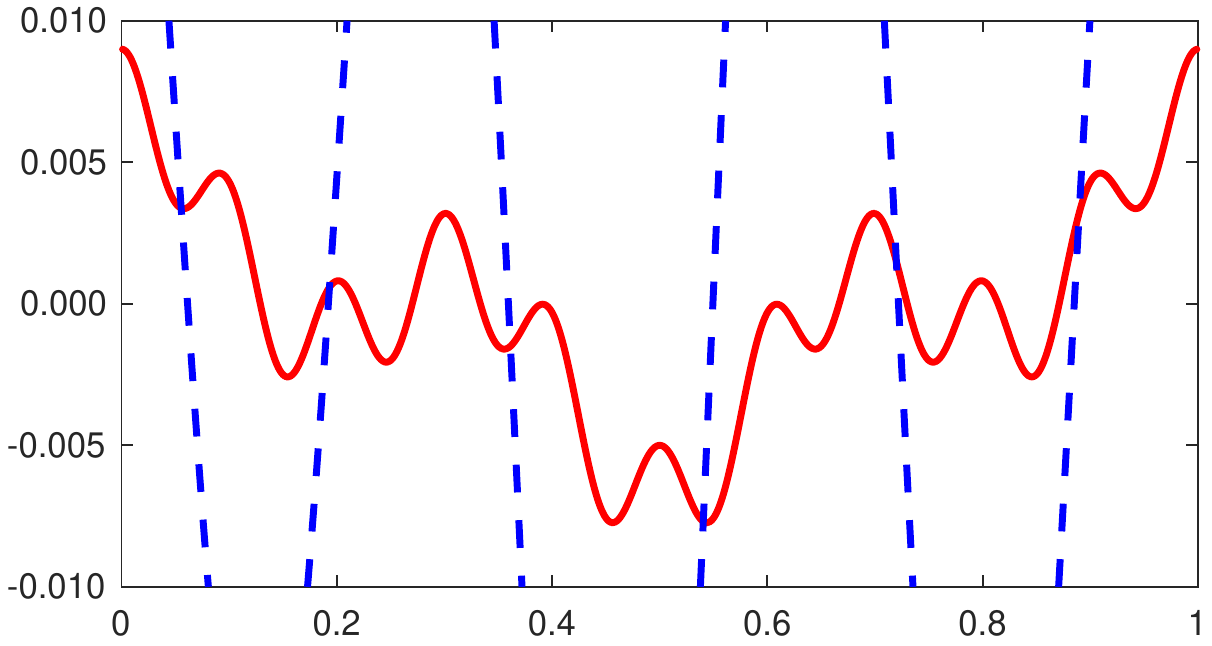} &
    \includegraphics[width=0.3\textwidth]{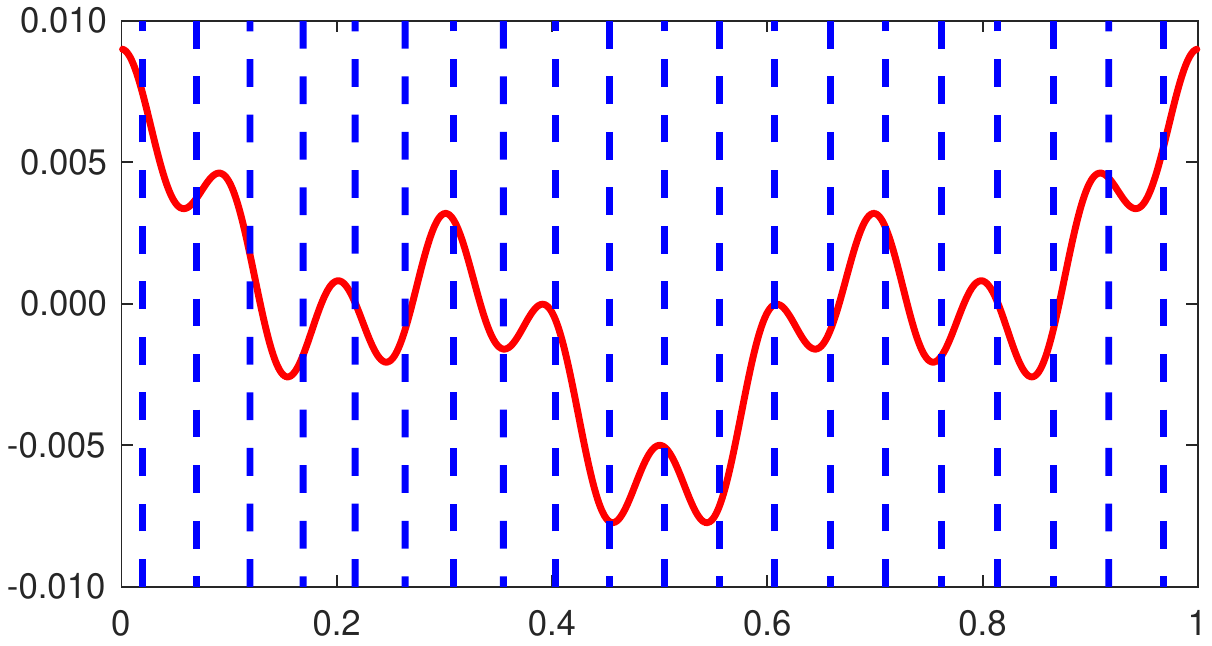} \\
    \includegraphics[width=0.3\textwidth]{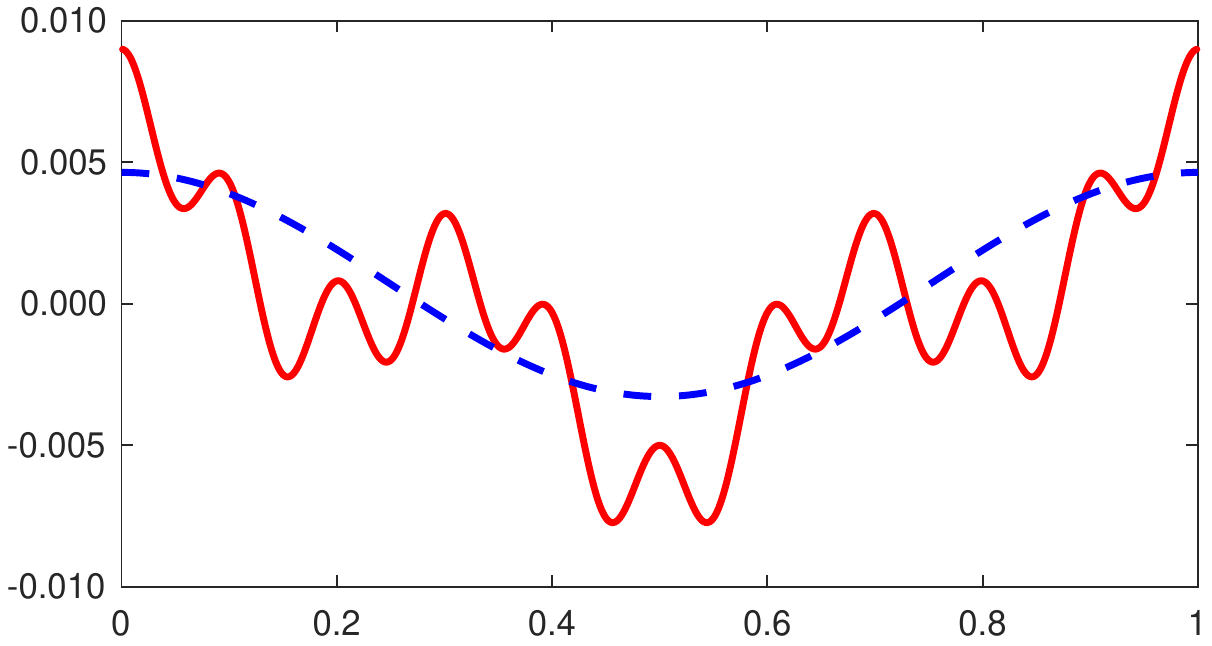} &
    \includegraphics[width=0.3\textwidth]{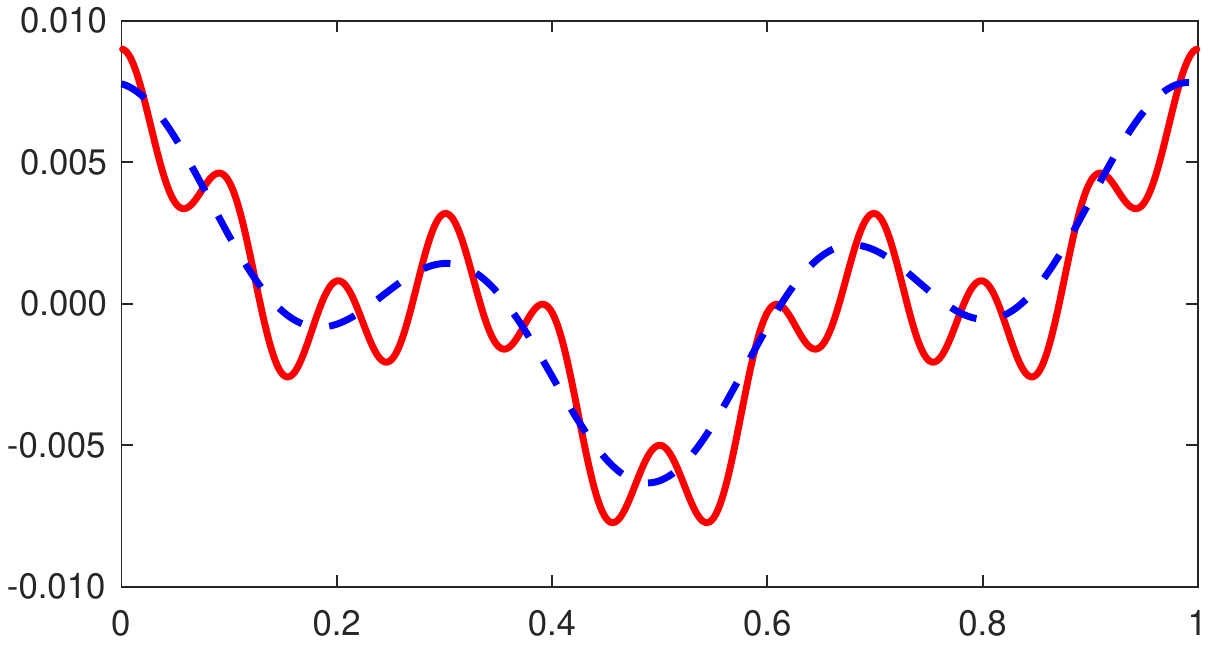} &
    \includegraphics[width=0.3\textwidth]{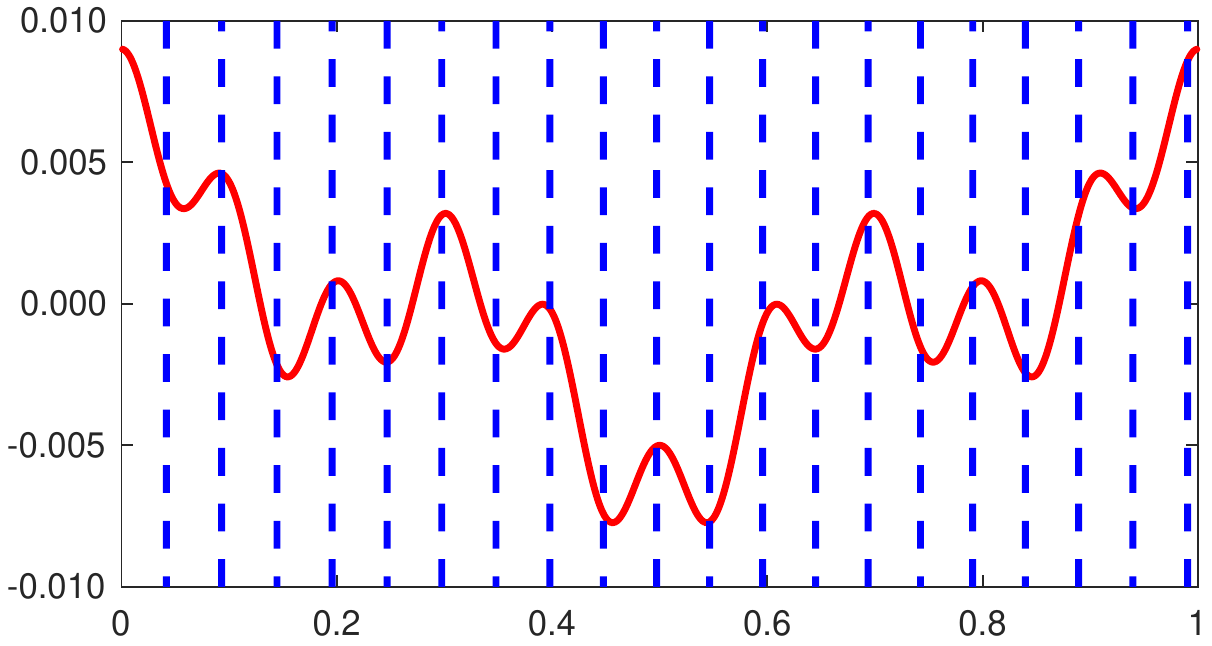} \\
    \includegraphics[width=0.3\textwidth]{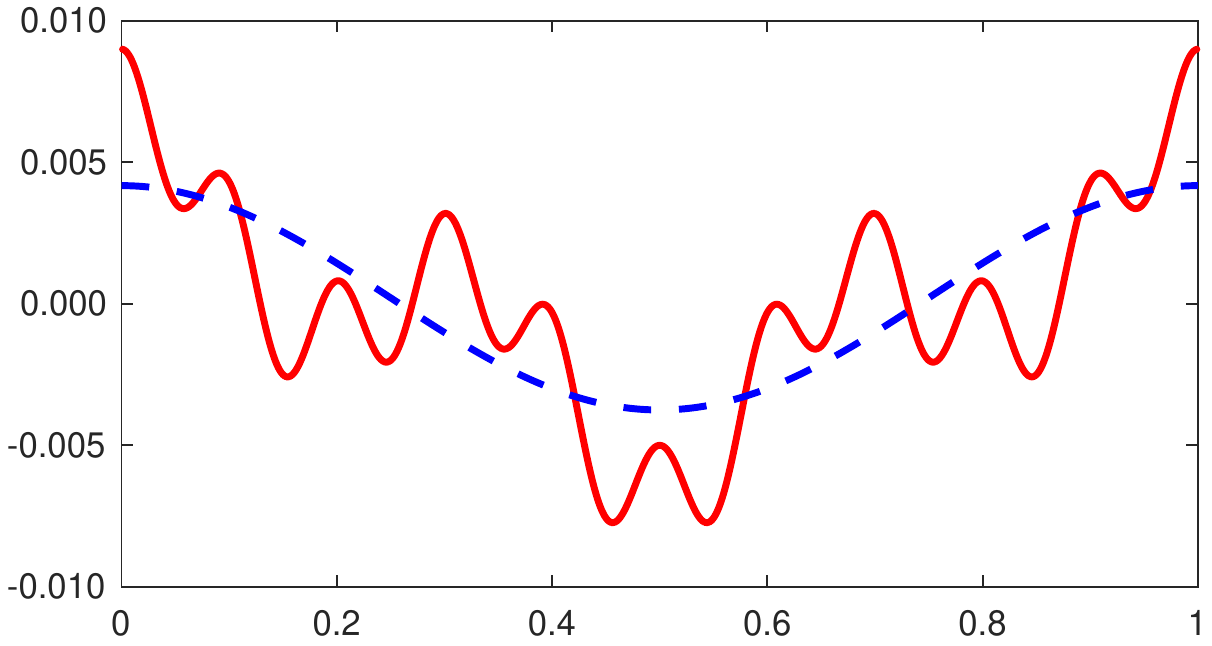} &
    \includegraphics[width=0.3\textwidth]{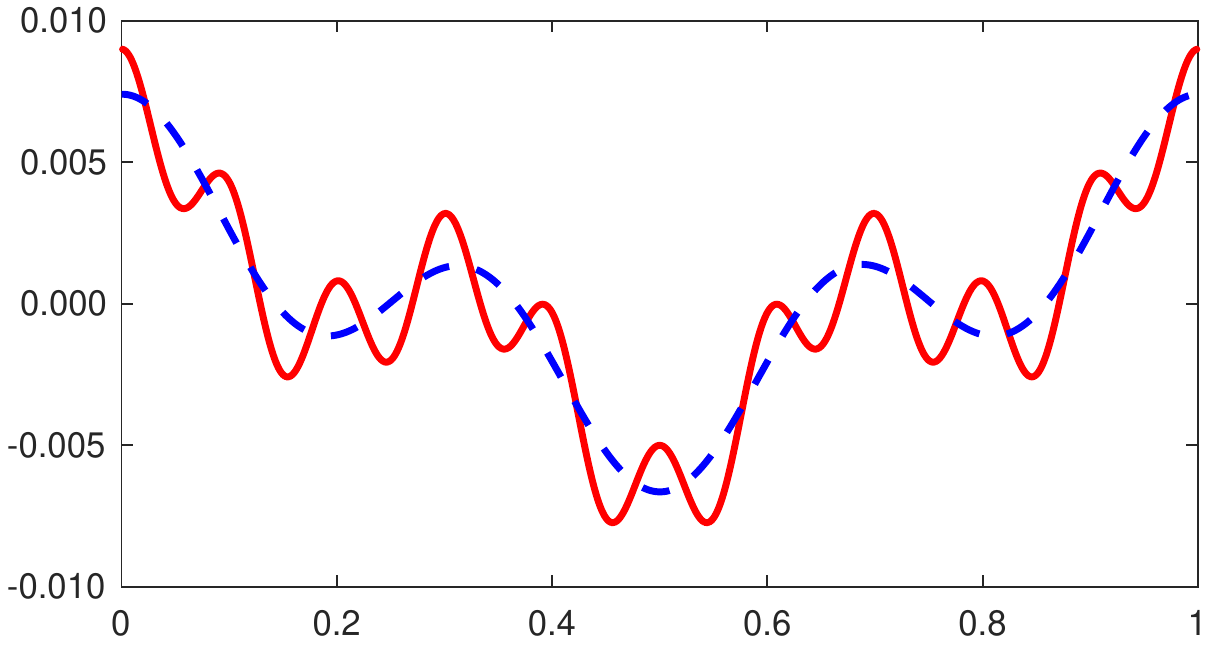} &
    \includegraphics[width=0.3\textwidth]{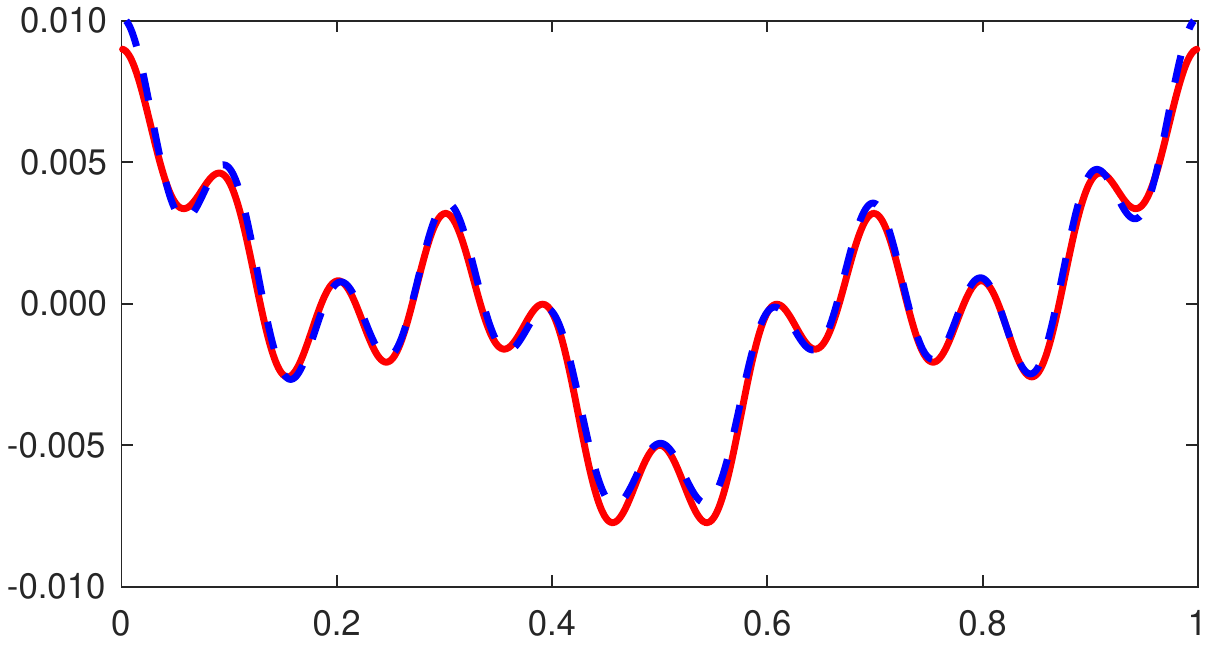} \\
    \includegraphics[width=0.3\textwidth]{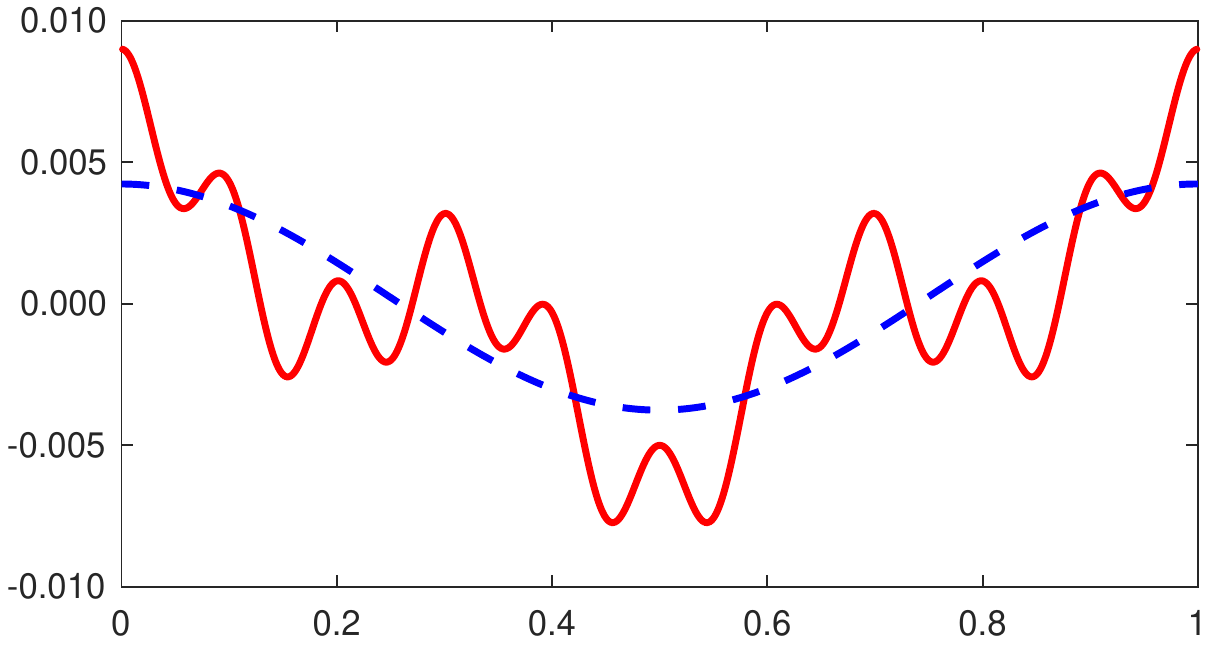} &
    \includegraphics[width=0.3\textwidth]{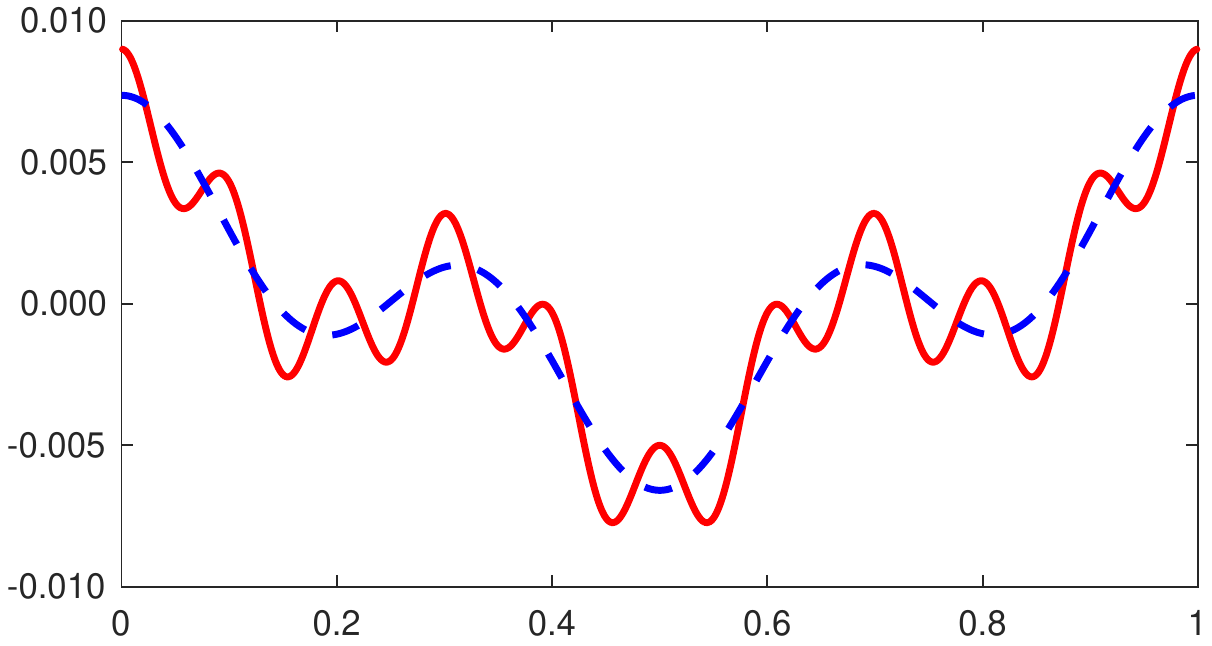} &
    \includegraphics[width=0.3\textwidth]{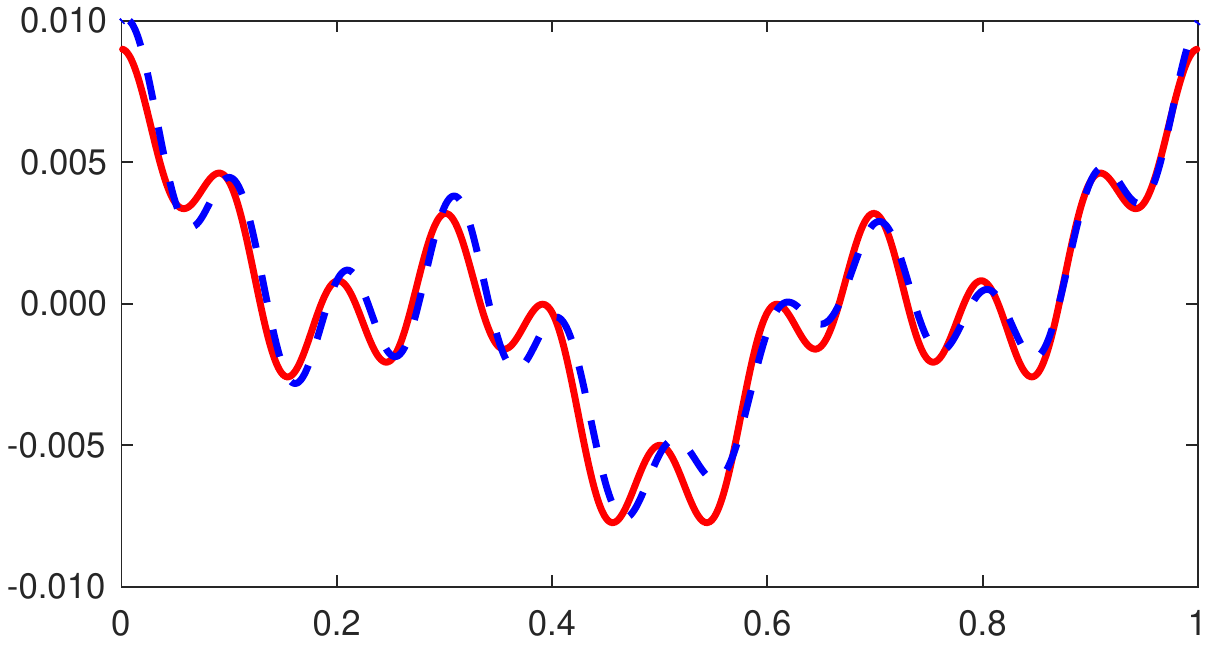} \\
    \includegraphics[width=0.3\textwidth]{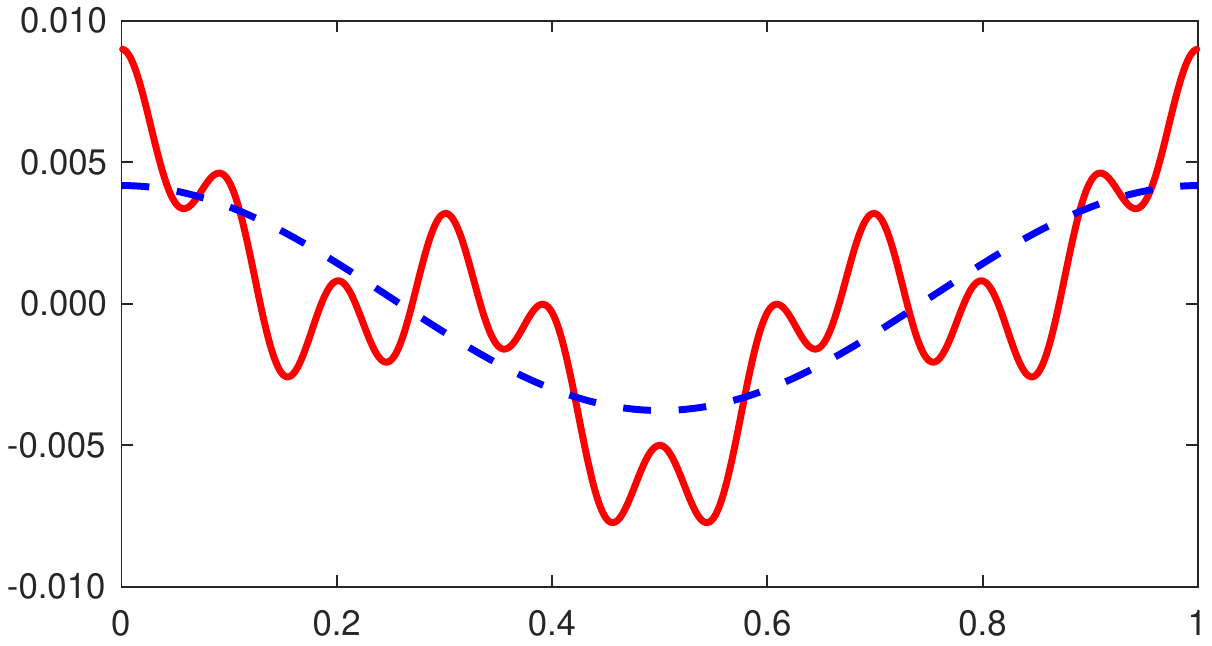} &
    \includegraphics[width=0.3\textwidth]{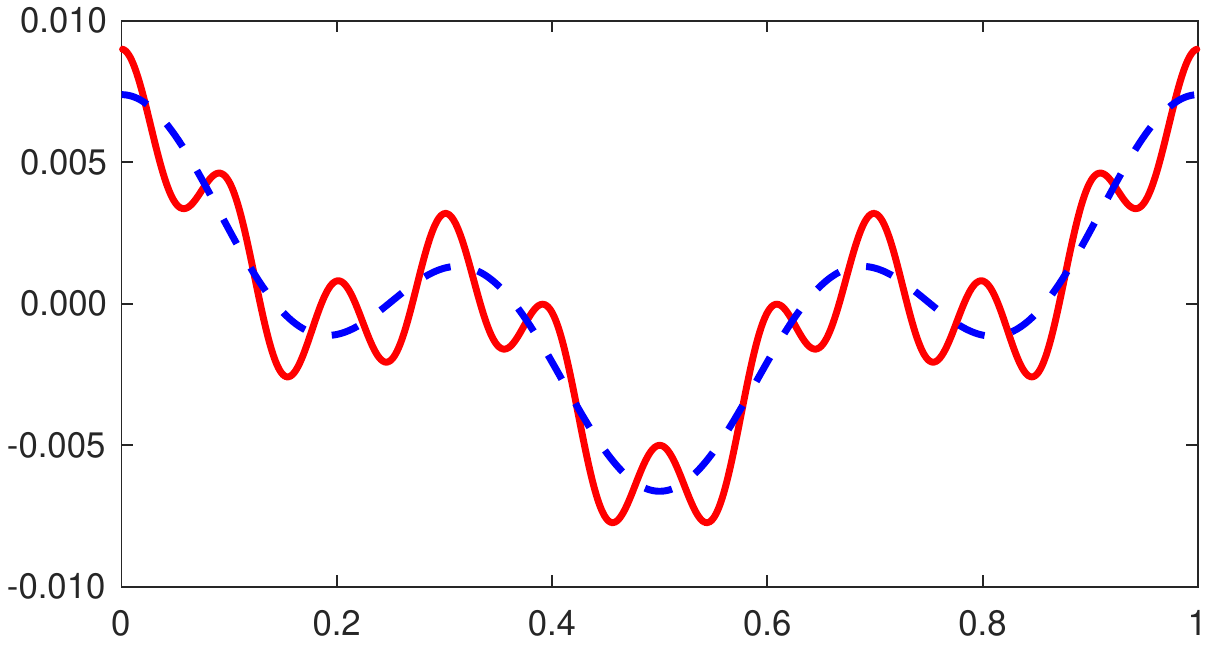} &
    \includegraphics[width=0.3\textwidth]{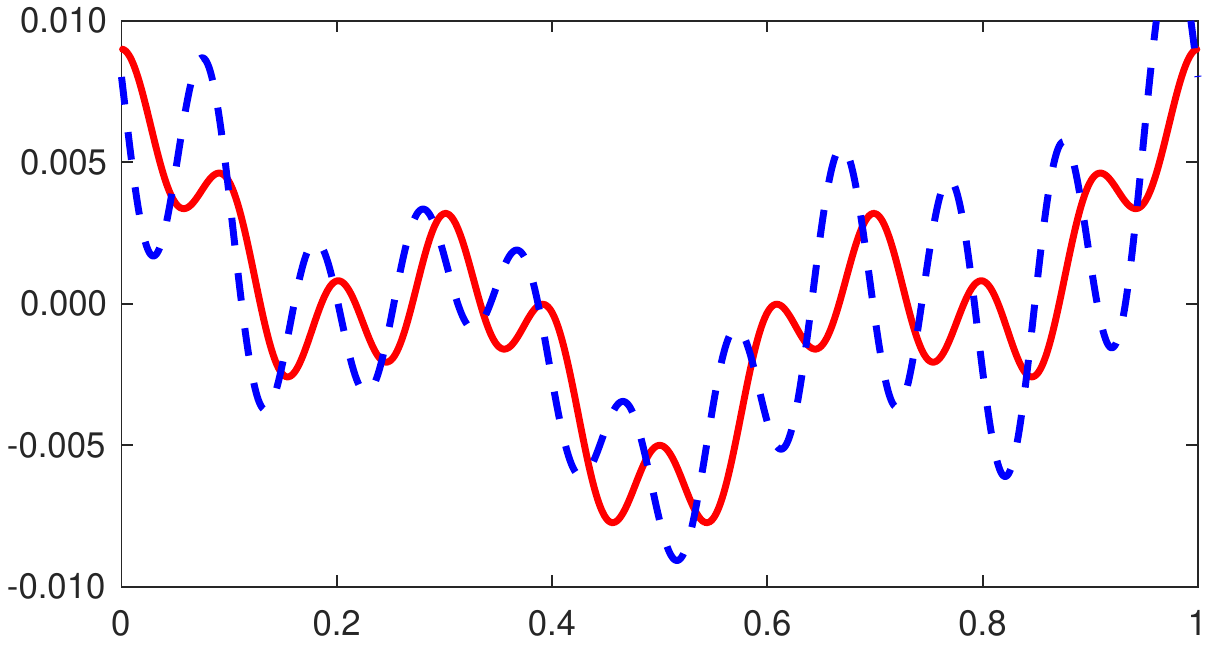}
  \end{tabular}
  \caption{Numerical experiments for the smooth profile function depicted in Eq. \eqref{eq:smoothProfile}. The red solid line represents the true profile, while the blue dashed line illustrates the reconstructed profile. The relative permittivity and permeability of the slab are specified as follows: Row 1: $\varepsilon = 1.0$ and $\mu = 1.0$; Row 2: $\varepsilon = 16.0$ and $\mu = 1.0$; Row 3: $\varepsilon = -1.0$ and $\mu = -1.0$; Row 4: $\varepsilon = -1 + 0.05i$ and $\mu = -0.97$; Row 5: $\varepsilon = -1 + 0.1i$ and $\mu = -1.06$. The cut-off frequency is specified as follows: Column 1: $N = 1$; Column 2: $N = 3$; Column 3: $N = 10$. \label{fig:smooth}}
\end{figure}

To recover the first, first two, or all three Fourier modes of the profile, the cut-off frequency $N$ must be set to at least $1$, $3$, or $10$, respectively. The reconstructions for the case without the slab are shown in Row 1 of Figure \ref{fig:smooth}. When $N=1$, only the first Fourier mode, $0.4 \cos (2 \pi t)$, is recovered. To recover the surface profile up to the 3rd Fourier mode, $N$ should be set to at least $3$. However, using $N=3$ leads to a solution with a significant error. Attempting to recover the 10th mode with $N=10$ proves unsuccessful. Better resolution can be achieved by decreasing the measurement distance $b$, but there is a lower limit for $b$, at least on top of the surface, and the error is not a monotonic function of $b$. A detailed error estimate in this case can be found in \cite{bao2014convergence}.

We can increase the resolution by introducing a dielectric slab with a high refractive index. Let us set $\varepsilon = 16$ and $\mu = 1$, which correspond to a refractive index of $4$. According to the analysis presented in \cite{bao2016nearFar}, the resolution can be increased by approximately a factor of $4$. Therefore, the profile can be reconstructed stably up to the 3rd Fourier mode, but not the 10th mode, as confirmed by the numerical results shown in Fig. \ref{fig:smooth} (Row 2).

We now replace the dielectric slab with a double negative metamaterial having perfectly matched parameters, i.e., $\varepsilon = \mu = -1$. According to the analysis presented in Section \ref{sec:perfect}, we should be able to recover all the frequency modes of the profile. This claim is verified by the numerical results shown in Figure \ref{fig:smooth} (Row 3), where the reconstructed profile remains unpolluted even when $N = 10$.

Lastly, we consider the case of imperfect parameters, taking $\varepsilon = - 1 + 0.05 i$ and $\mu = - 0.97$ as an example. The corresponding reconstructions are shown in Fig. \ref{fig:smooth} (Row 4). As expected, the results deteriorate with increasing deviation from the ideal case, particularly for $N = 10$, but they are still relatively close to the ideal case. The reconstruction accuracy will naturally decrease as $\varepsilon$ and $\mu$ move further away from $-1$. Fig. \ref{fig:smooth} (Row 5) displays the reconstruction with $\varepsilon = - 1 + 0.1 i$ and $\mu = - 1.06$. Note that the first and third modes are still well recovered, while the 10th mode is not.

\begin{figure}[t]
\centering
\includegraphics[width=0.55\textwidth]{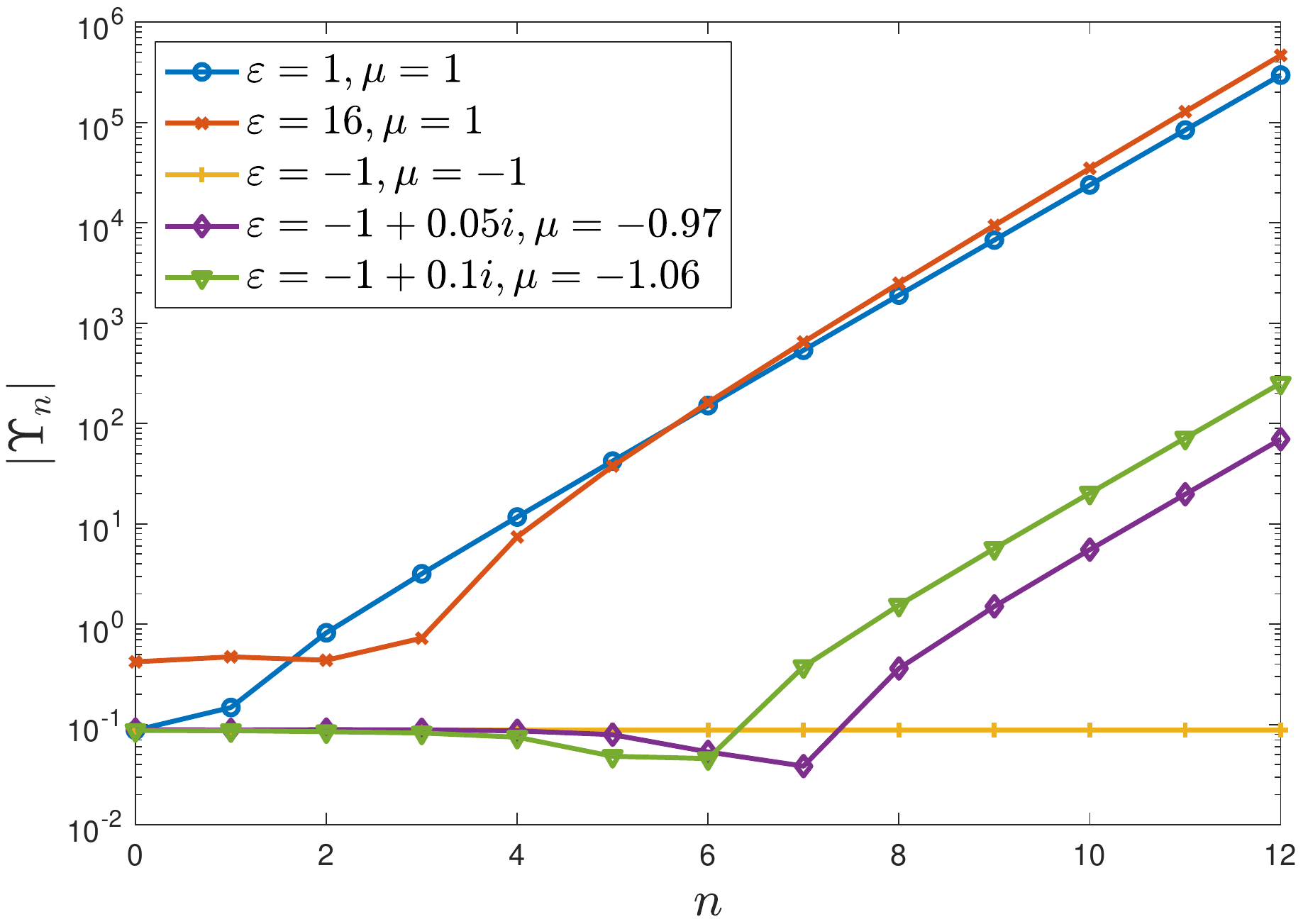}
\caption{The scaling factor $|\Upsilon_n|$ in Eq. \eqref{eq:scaling} is plotted on a logarithmic scale for various values of $\varepsilon$ and $\mu$ in the numerical experiments presented in Fig. \ref{fig:smooth}. \label{fig:Upsilon}}
\end{figure}

In Fig. \ref{fig:Upsilon}, we present the absolute value of the scaling factor $|\Upsilon_n|$ defined in Eq. \ref{eq:scaling} plotted against the mode number $n$ for the various values of $\varepsilon$ and $\mu$ used in the previous numerical experiments. We make the following observations:
\begin{enumerate}[(i)]
\item  No slab ($\varepsilon = 1, \mu = 1$): $|\Upsilon_n|$ increases exponentially from $n=0$.
\item High refractive index ($\varepsilon = 16, \mu = 1$): $|\Upsilon_n|$ remains small for a few small values of $n$ before increasing exponentially with $n$, resulting in enhanced resolution compared to case (i).
\item Perfectly matched parameters ($\varepsilon = \mu = - 1$): $|\Upsilon_n|$ remains small for all $n$, leading to unlimited resolution.
\item Imperfect parameters ($\varepsilon = - 1 + 0.05 i, \mu = - 0.97$ and $\varepsilon = - 1 + 0.1 i, \mu = - 1.06$): Similar to case (ii), $|\Upsilon_n|$ remains small for a few small values of $n$ before increasing exponentially. Hence the resolution is enhanced, and it would continue to improve as $\varepsilon \rightarrow -1$ and $\mu \rightarrow -1$.
\end{enumerate}

\subsection{Nonsmooth Profiles}

Although the derivation of the reconstruction formula requires the existence of $g'$ and $g''$, our method can still be applied to nonsmooth profiles due to the approximation properties of Fourier series. We first consider a nonsmooth profile defined by \(f(x) = \delta g(x)\), where \(\delta = 0.01\) and
\begin{equation} \label{eq:nonsmooth}
  g(x) = (1 - 10 | x - 0.3 |) \chi_{[0.2, 0.4]}(x) + (1 - 10 | x - 0.7 |)
  \chi_{[0.6, 0.8]}(x).
\end{equation}

\begin{figure}[t]
  \begin{tabular}{ccc}
    \includegraphics[width=0.3\textwidth]{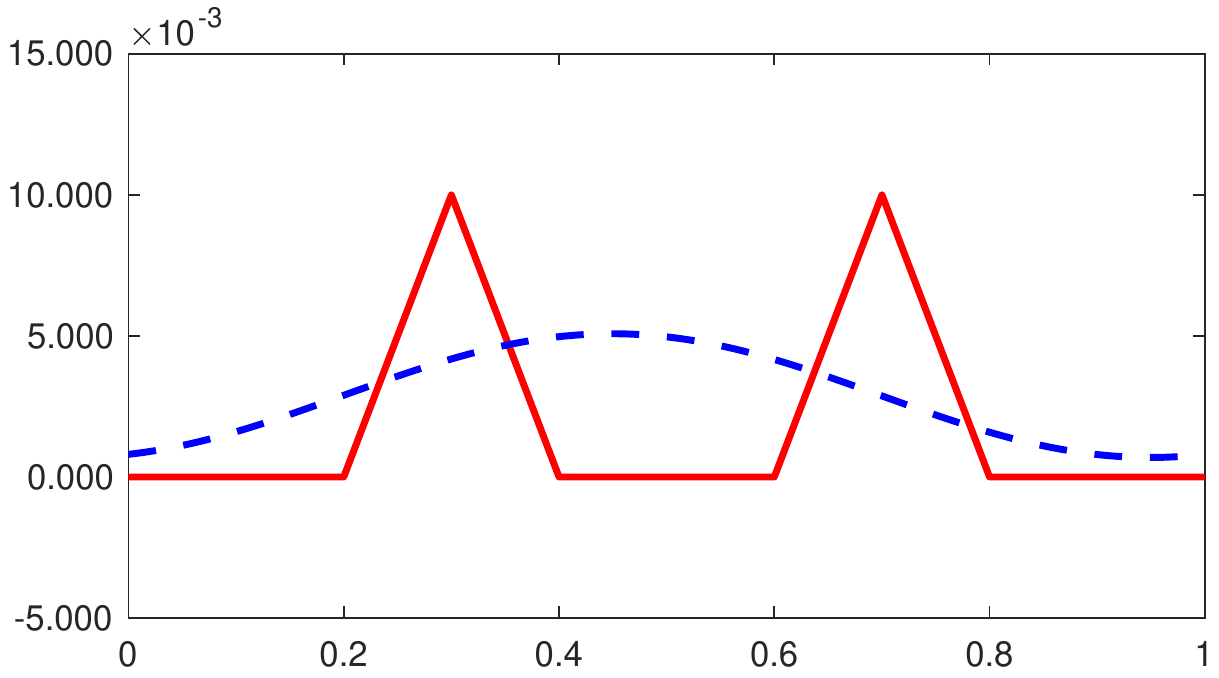} &
    \includegraphics[width=0.3\textwidth]{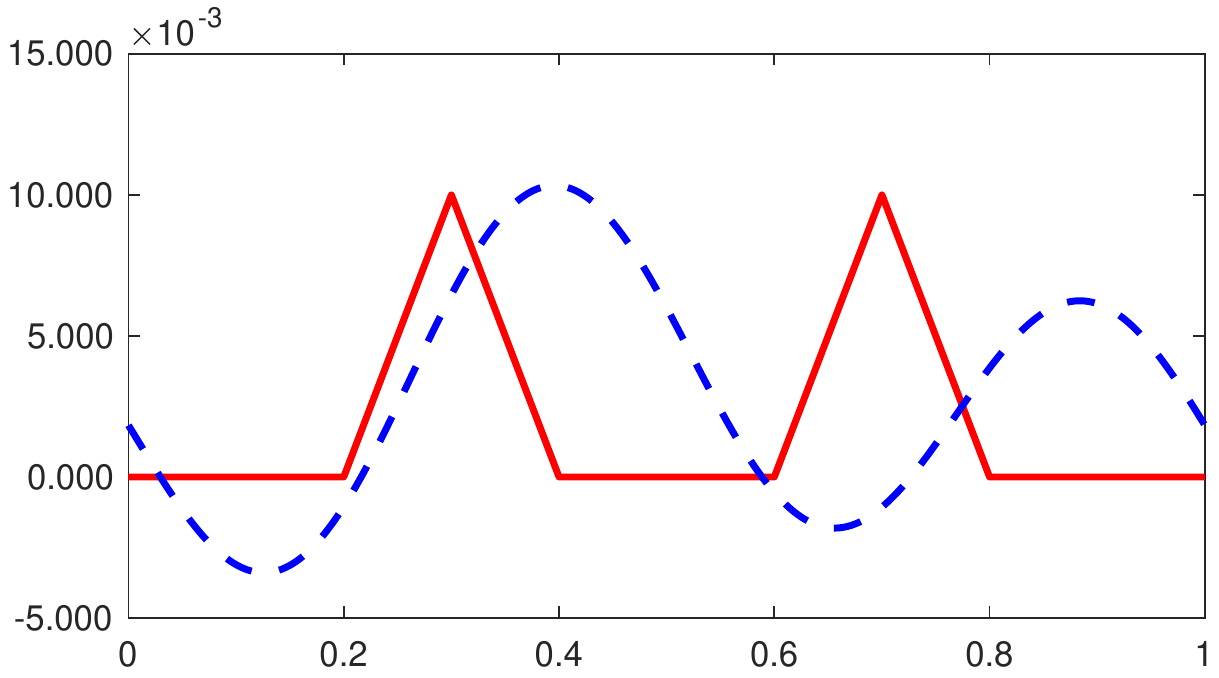} &
    \includegraphics[width=0.3\textwidth]{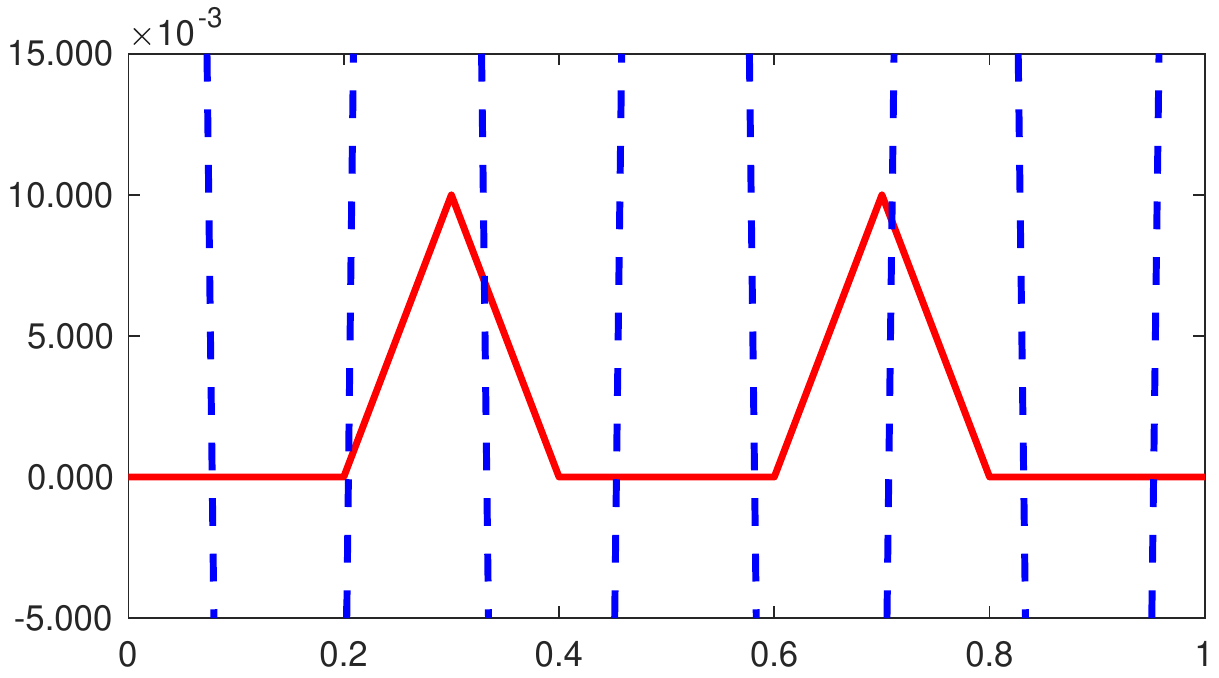} \\
    \includegraphics[width=0.3\textwidth]{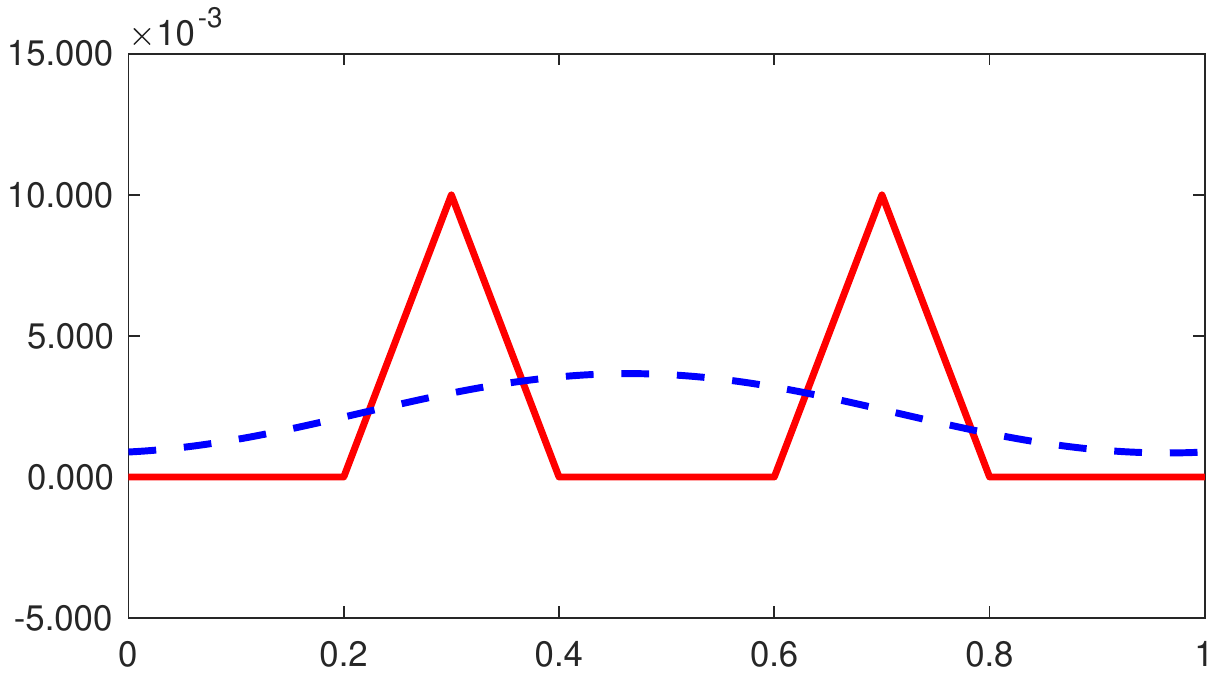} &
    \includegraphics[width=0.3\textwidth]{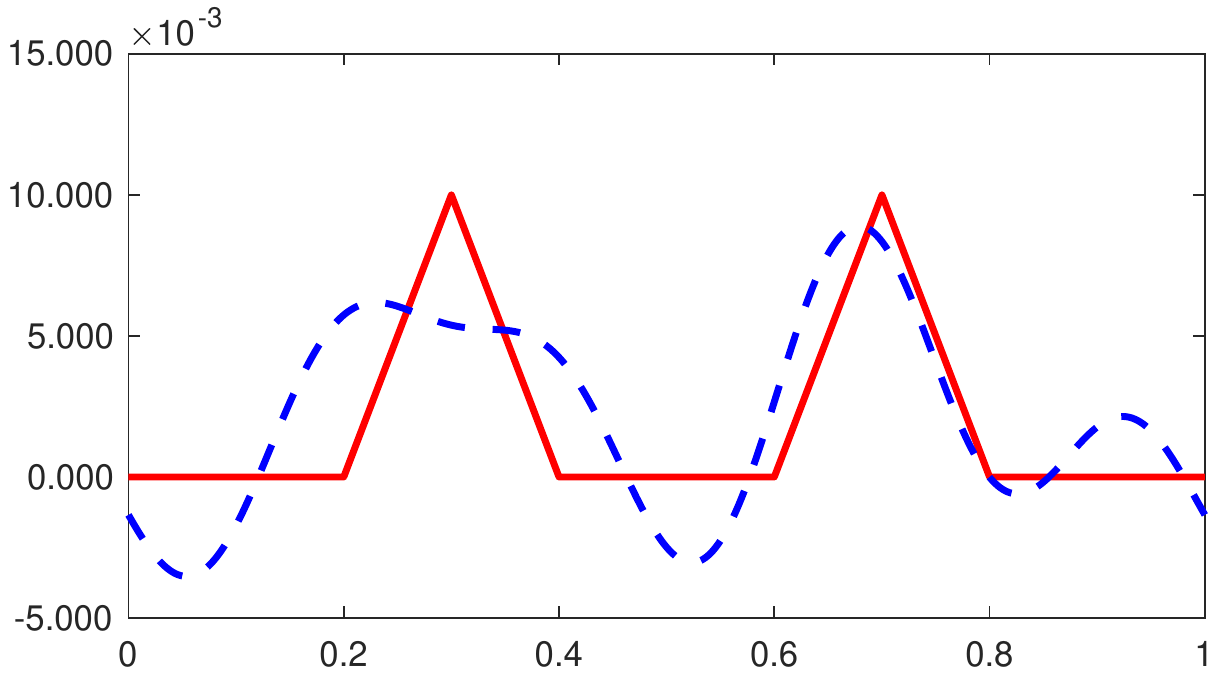} &
    \includegraphics[width=0.3\textwidth]{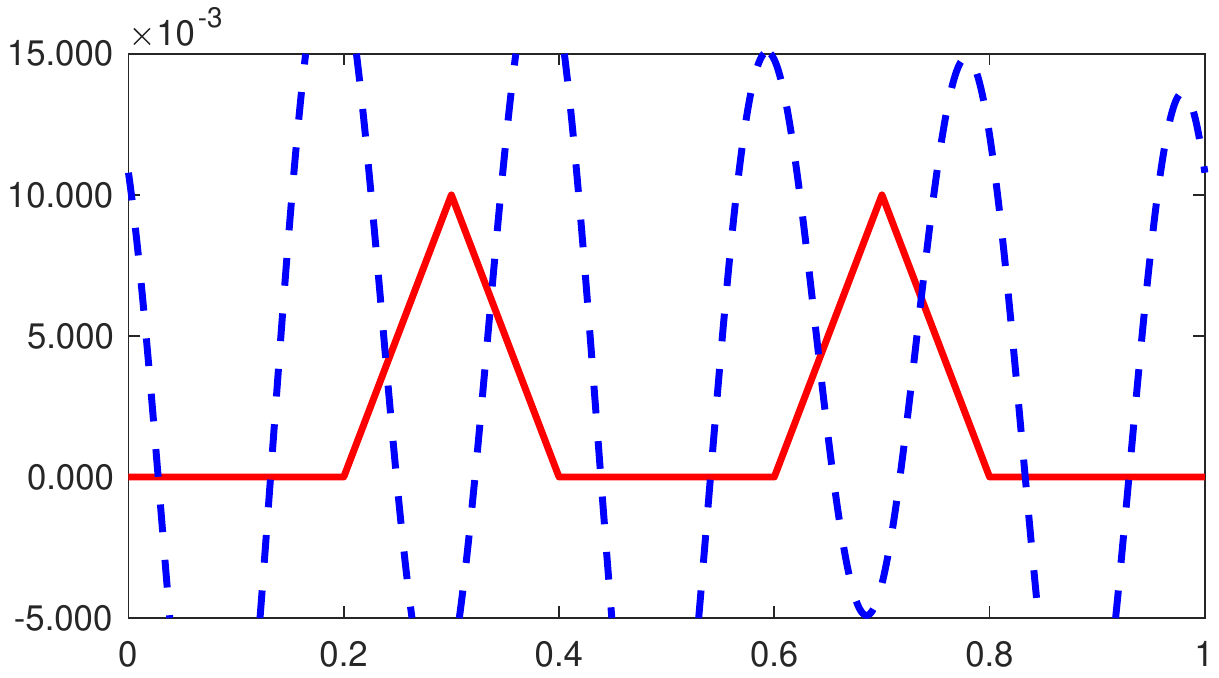} \\
    \includegraphics[width=0.3\textwidth]{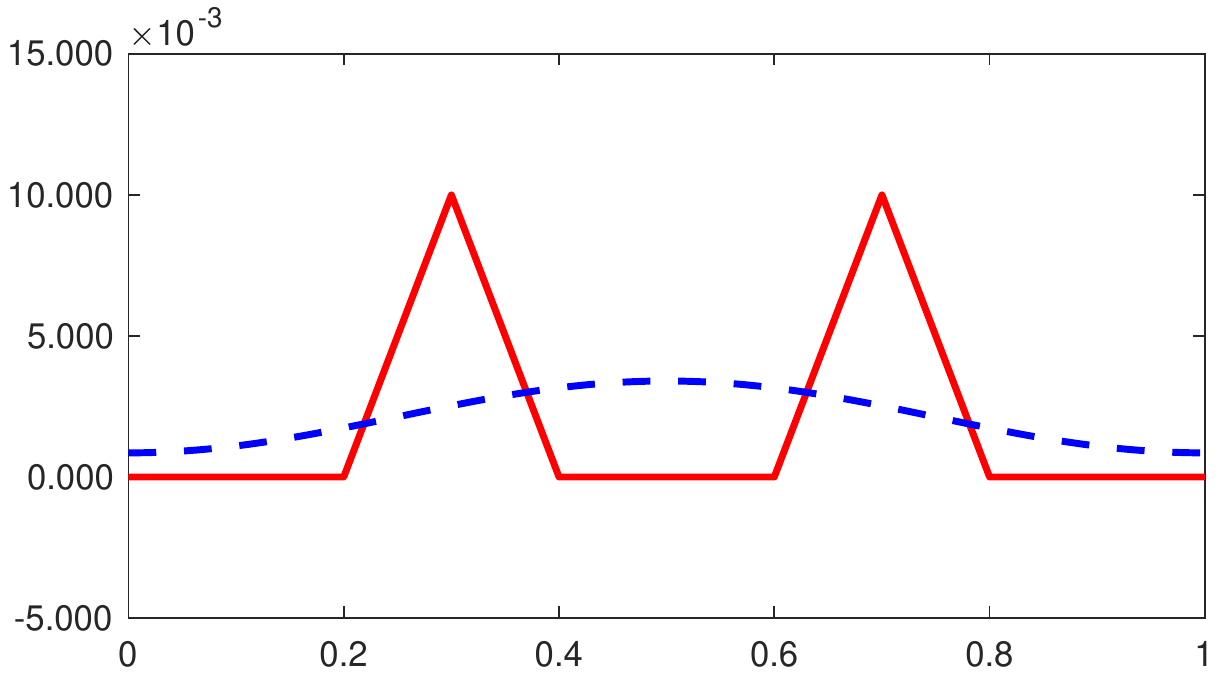} &
    \includegraphics[width=0.3\textwidth]{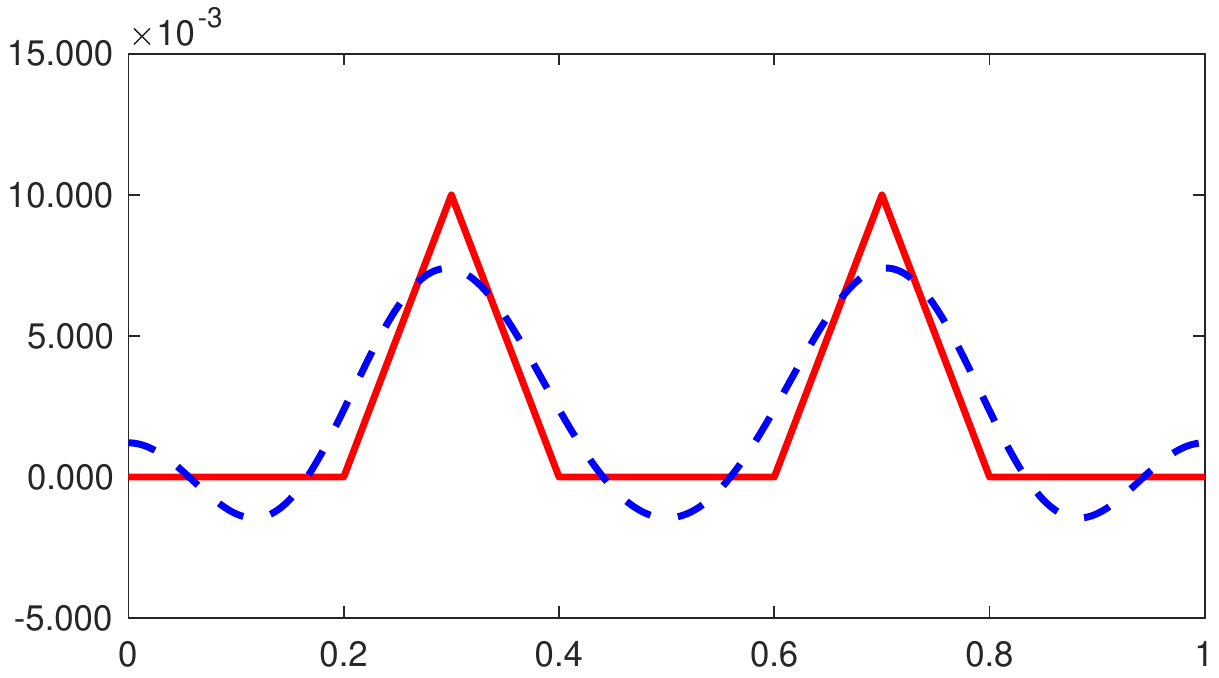} &
    \includegraphics[width=0.3\textwidth]{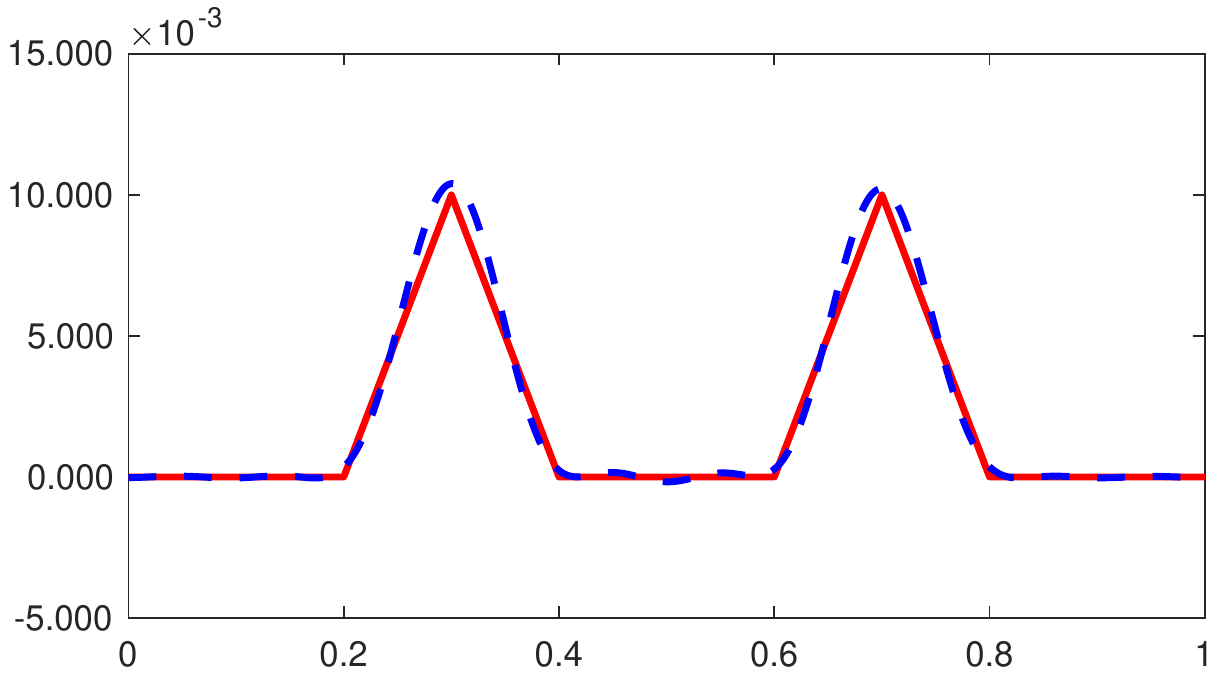}
  \end{tabular}
  \caption{Numerical experiments for the nonsmooth profile function \eqref{eq:nonsmooth}. The red solid line represents the true profile, while the blue dashed line corresponds to the reconstructed profile. Row 1: \(\varepsilon = 1.0\), \(\mu = 1.0\), \(N = 1, 2, 4\); Row 2: \(\varepsilon = 16.0\), \(\mu = 1.0\), \(N = 1, 4, 5\); Row 3: \(\varepsilon = -1.0\), \(\mu = -1.0\), \(N = 1, 4, 8\). \label{fig:nonsmooth}}
\end{figure}

Note that this profile contains infinitely many nonzero Fourier modes. The results of the experiment are displayed in Fig. \ref{fig:nonsmooth}. For the case without a slab (\(\varepsilon = 1.0\), \(\mu = 1.0\)), we begin with the cut-off frequency \(N = 1\). The reconstructed profile is smooth, but the two humps are indistinguishable. With \(N = 2\), the reconstruction starts capturing the humps, but the accuracy remains low. When we attempt to increase the accuracy by raising \(N\) to \(4\), the reconstruction becomes significantly contaminated by noise. 

Next, we insert a slab with a high refractive index (\(\varepsilon = 16.0\), \(\mu = 1.0\)). The reconstruction now becomes stable and more accurate with \(N = 4\). However, increasing \(N\) to \(5\) leads to unstable results again. 

Finally, we replace the slab with a superlens (\(\varepsilon = - 1.0\), \(\mu = - 1.0\)). The reconstruction is now stable and accurate up to at least \(N = 8\).

\begin{figure}[t]
  \begin{tabular}{ccc}
    \includegraphics[width=0.3\textwidth]{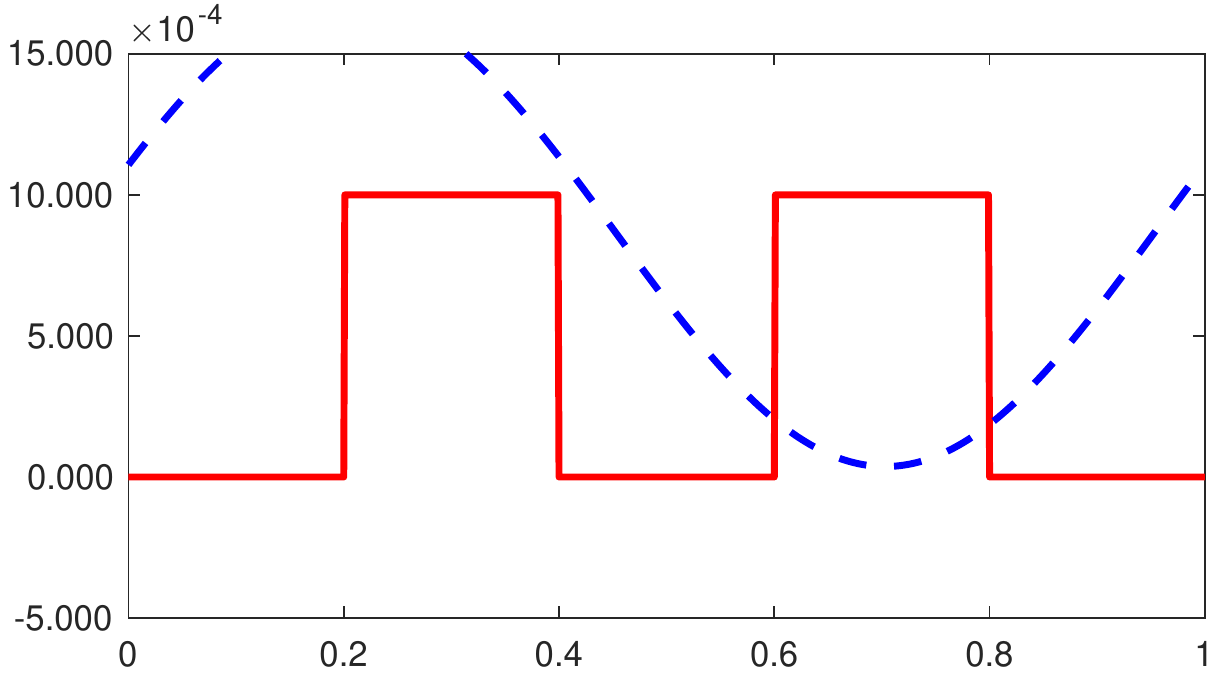} &
    \includegraphics[width=0.3\textwidth]{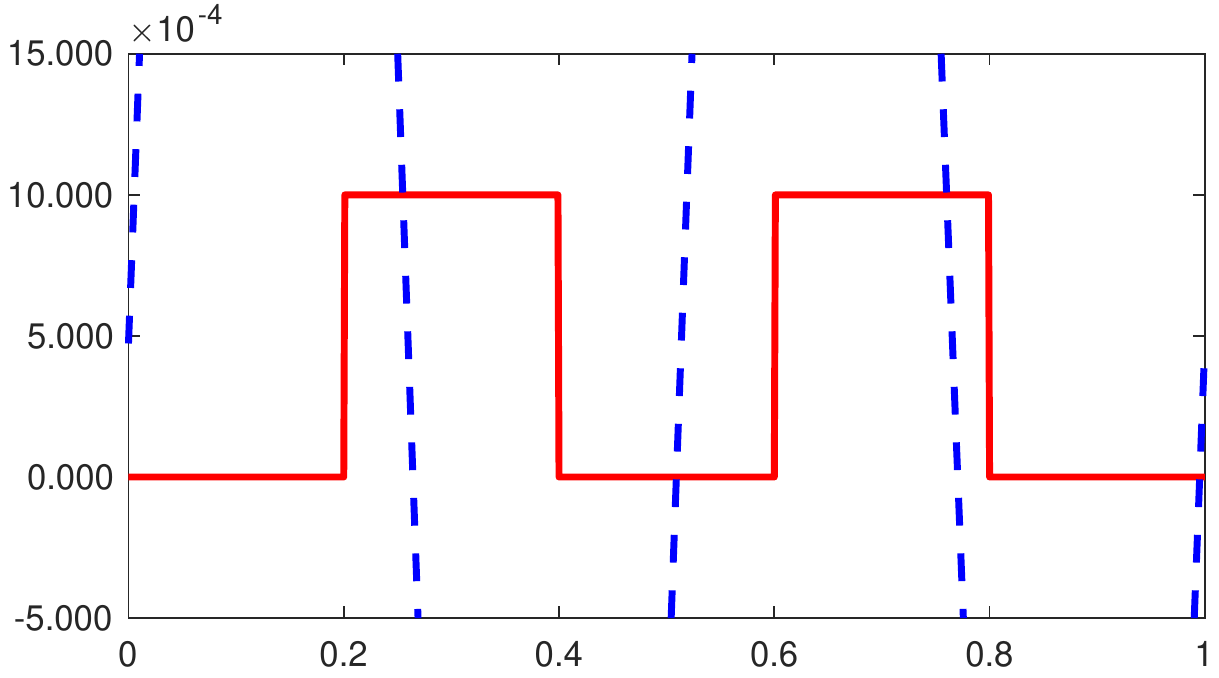} &
    \includegraphics[width=0.3\textwidth]{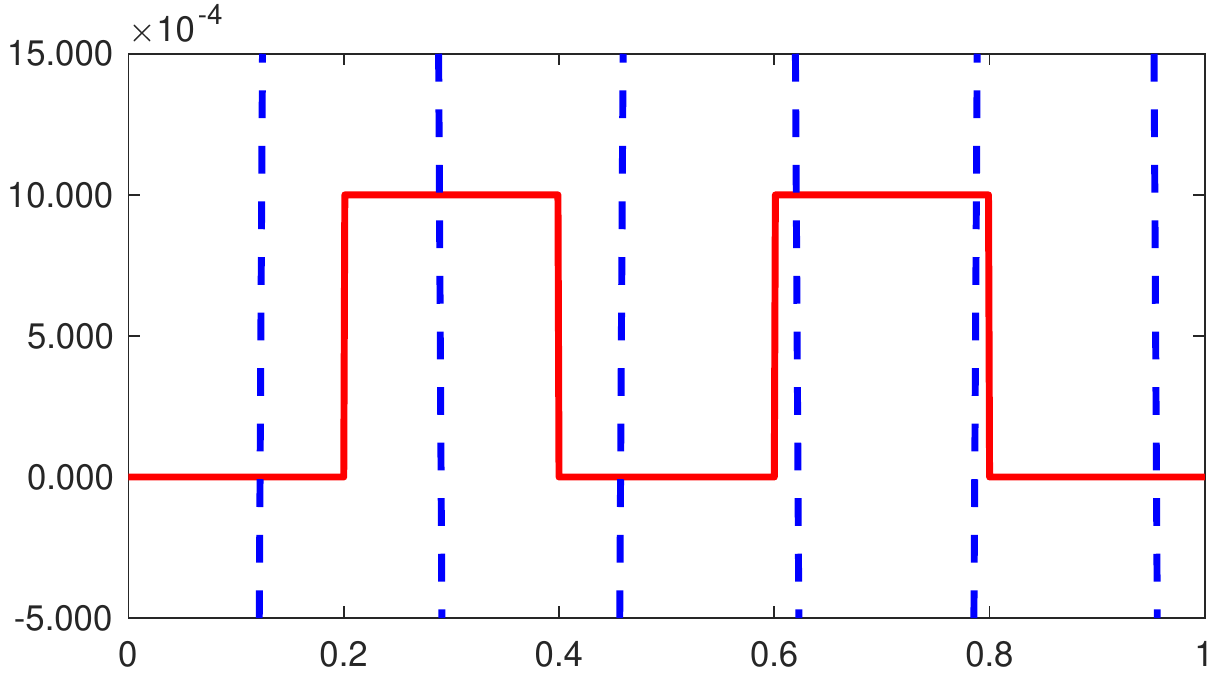} \\
    \includegraphics[width=0.3\textwidth]{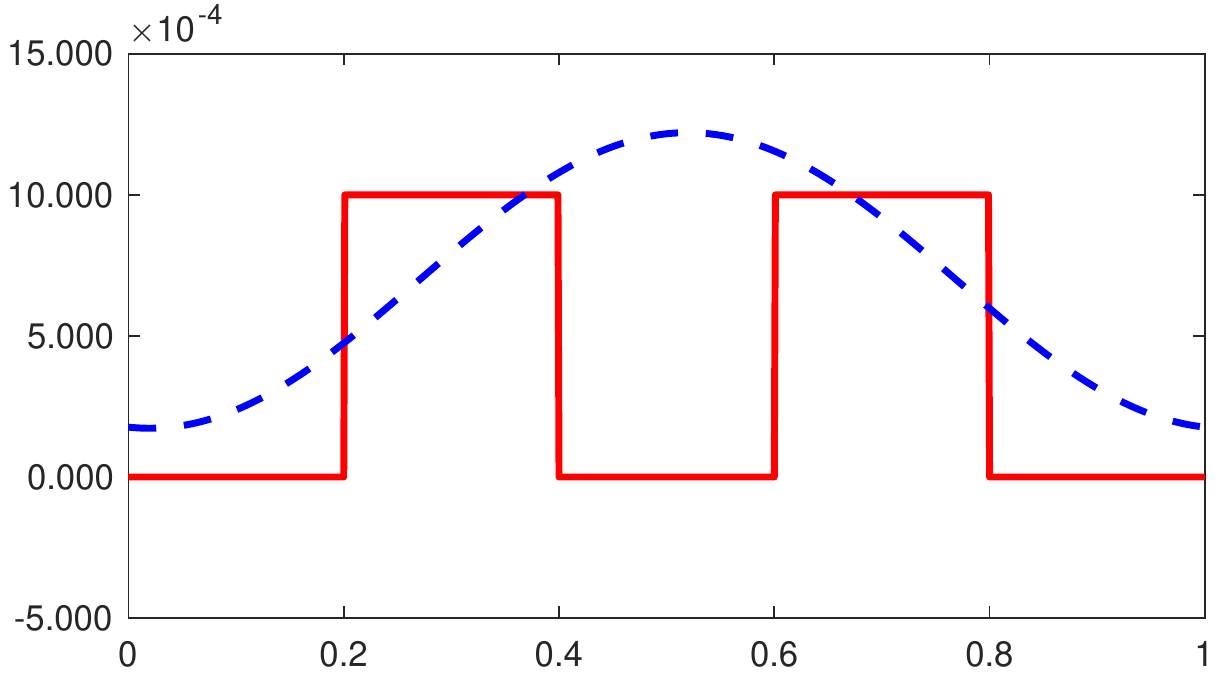} &
    \includegraphics[width=0.3\textwidth]{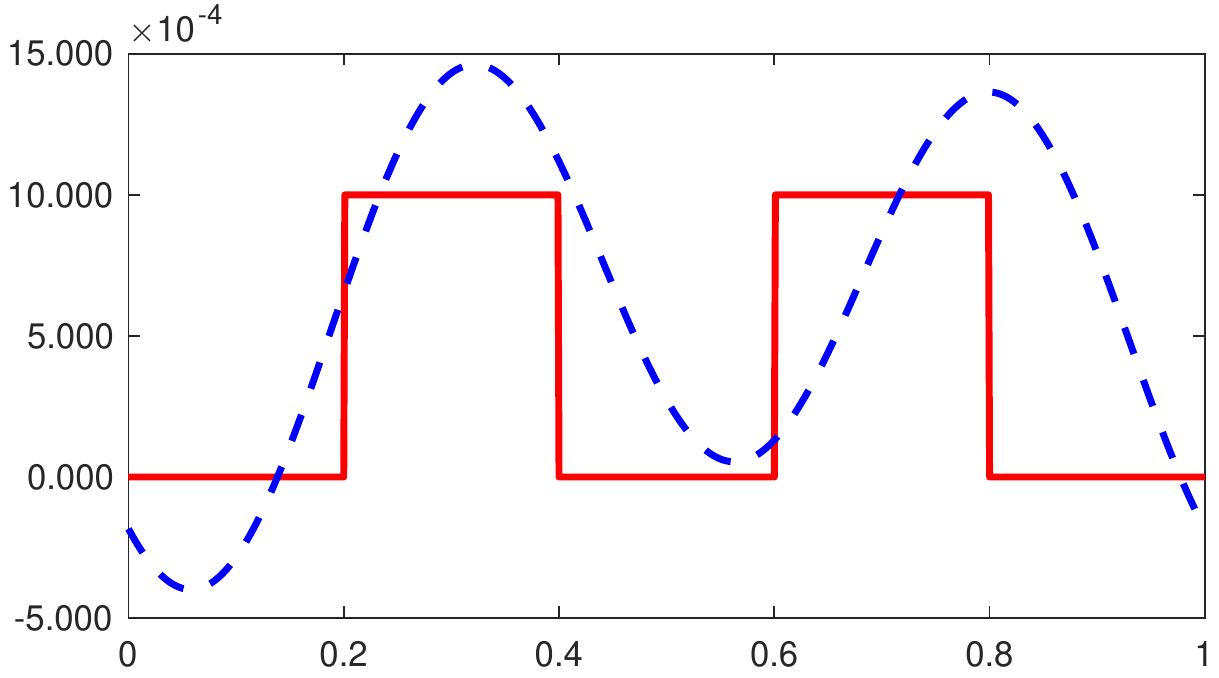} &
    \includegraphics[width=0.3\textwidth]{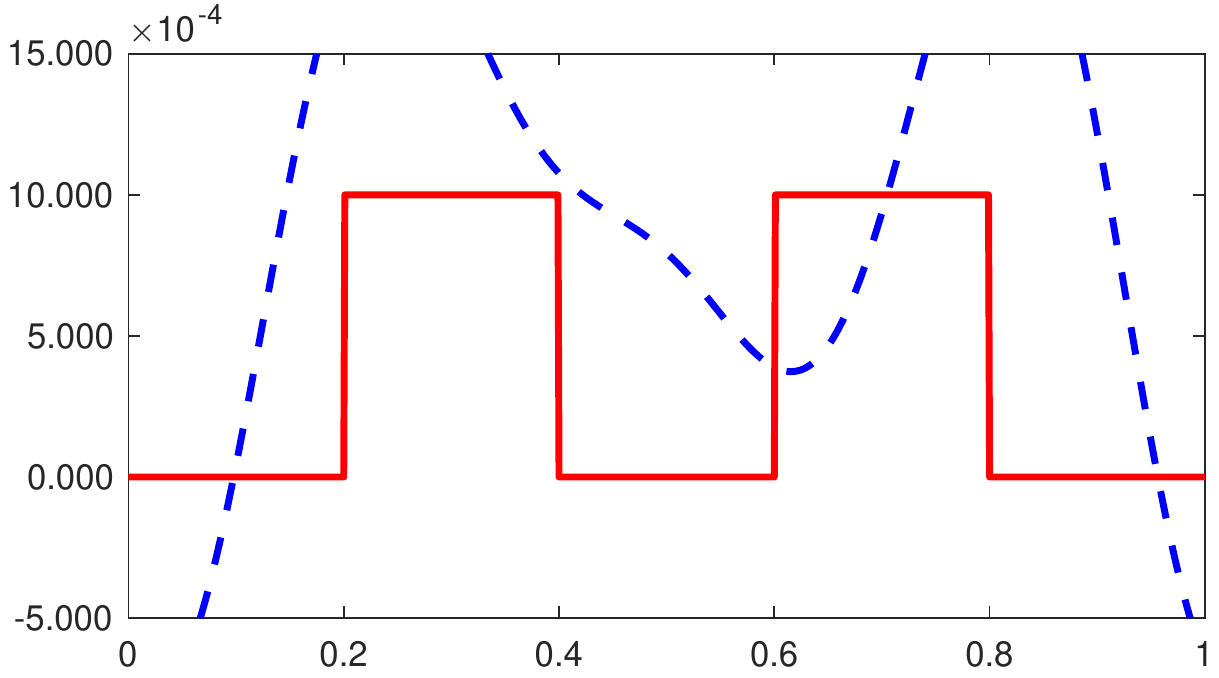} \\
    \includegraphics[width=0.3\textwidth]{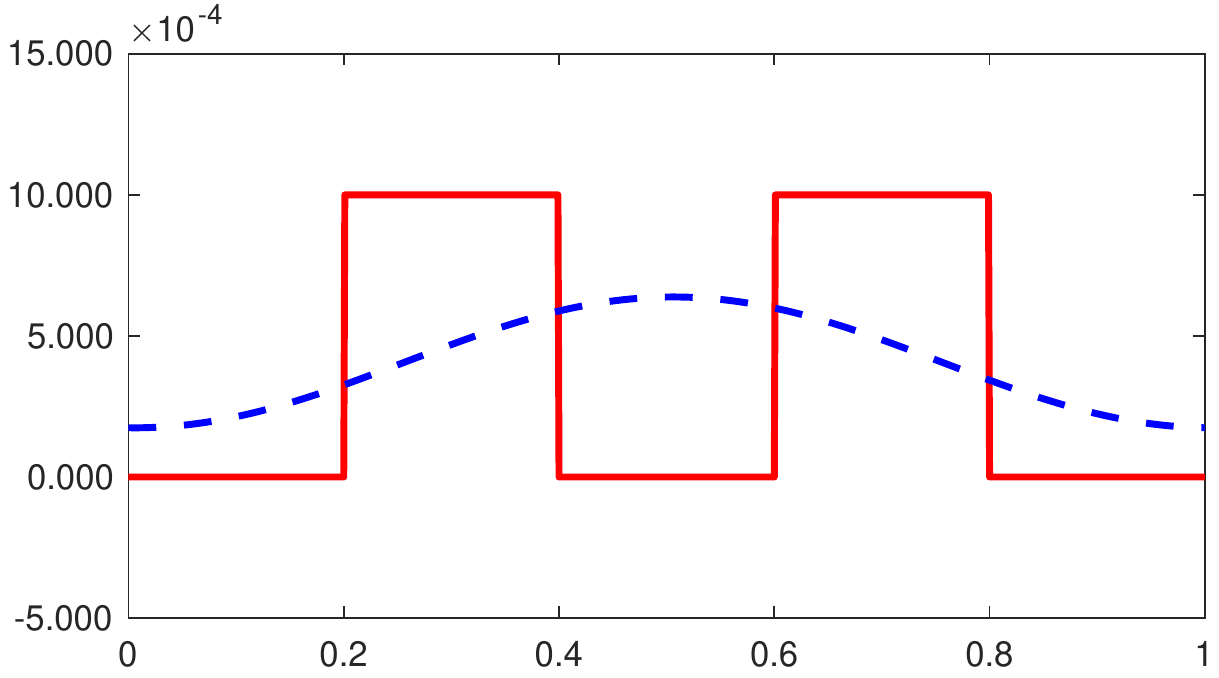} &
    \includegraphics[width=0.3\textwidth]{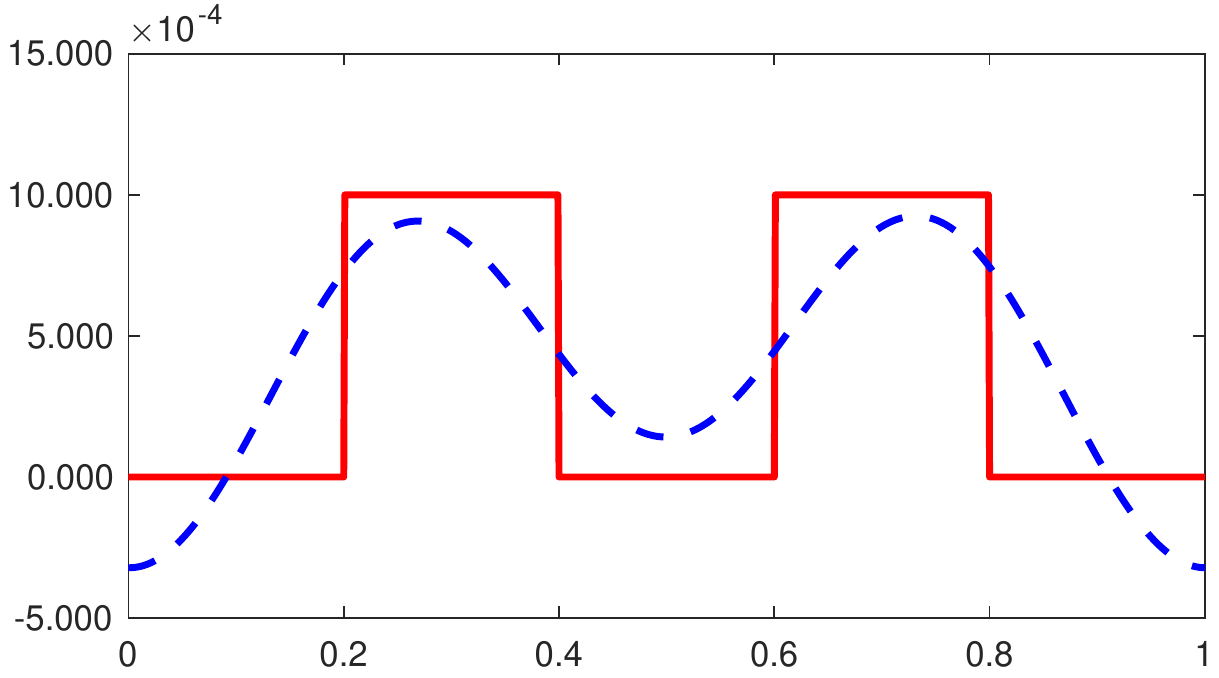} &
    \includegraphics[width=0.3\textwidth]{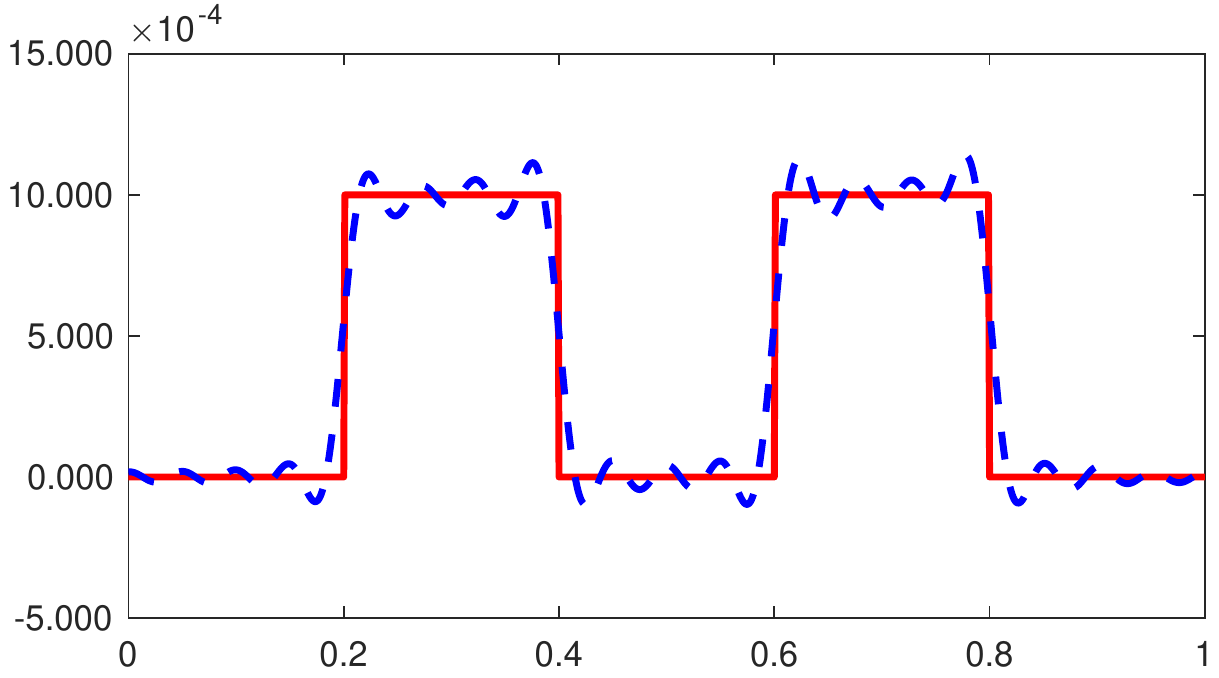}
  \end{tabular}
  \caption{Numerical experiments for the nonsmooth profile function \eqref{eq:discontinuous}. The red solid line represents the true profile, while the blue dashed line corresponds to the reconstructed profile. Row 1: $\varepsilon = 1.0, \mu = 1.0, N = 1, 2, 3$;
  Row 2: $\varepsilon = 16.0, \mu = 1.0, N = 1, 2, 3$; Row 3: $\varepsilon = -
  1.0, \mu = - 1.0, N = 1, 2, 20$. \label{fig:discontinuous}}
\end{figure}

Next, we consider a discontinuous profile given by $f (x) = \delta g (x)$,
where $\delta = 0.001$ and
\begin{equation}
  g (x) = \chi_{[0.2, 0.4]} (x) + \chi_{[0.6, 0.8]} (x) .
  \label{eq:discontinuous}
\end{equation}
We conducted another experiment using a setup similar to the previous one and present the results in Fig. \ref{fig:discontinuous}. In the case without a slab, with \(\varepsilon = 1.0\) and \(\mu = 1.0\), we did not obtain satisfactory results for any cut-off frequency, \(N\). For the slab with a high refractive index characterized by \(\varepsilon = 16.0\) and \(\mu = 1.0\), the optimal reconstruction is achieved with \(N = 2\), while the quality deteriorates for \(N \geqslant 3\). In the superlens case with \(\varepsilon = -1.0\) and \(\mu = -1.0\), the reconstruction remains stable up to at least \(N = 20\). However, the Gibbs phenomenon becomes more pronounced for larger values of \(N\).

\section{Conclusion}\label{sec:conclusion}

In this paper, we investigated the superresolution effect of double negative index materials using an inverse scattering problem approach for diffraction gratings with low amplitude profiles. By applying the transformed field expansion method, we derived a simple reconstruction formula for the linear approximation of the scattering problem. Our findings suggest that it is feasible to achieve unlimited resolution with perfectly matched values of permittivity and permeability. To demonstrate the effectiveness of our reconstruction method, we conducted numerical experiments for both smooth and nonsmooth profiles with perfect or imperfect parameters.

There are several directions for future research on this topic. Firstly, the direct scattering problem with sign-changing coefficients for periodic structures requires further exploration to establish its well-posedness. The uniqueness of the inverse problem can also be investigated, and the convergence and error estimates of our proposed method can be analyzed. Moreover, it is possible to extend our method to other forms of waves and to the three-dimensional scenarios. For more general profile or inverse medium scattering problems, optimization or sampling methods may prove effective. Finally, to apply our method in real-world scenarios, it is crucial to incorporate the metamaterial's constitution and measuring device into the model.

\section*{Acknowledgments}

The research of PL was supported in part by the NSF grant DMS-2208256.

\bibliographystyle{unsrt}
\bibliography{superlens}

\end{document}